\documentclass[12pt,reqno,final]{amsart}
\usepackage[backref,colorlinks,linkcolor=red,anchorcolor=green,citecolor=blue]{hyperref}
\usepackage{mathrsfs}
\usepackage{amssymb}
\usepackage{cases}
\usepackage{amsmath}
\usepackage{stmaryrd}
\usepackage{txfonts}
\usepackage{indentfirst,latexsym,bm}
\usepackage{graphicx}
\usepackage{psfrag}
\usepackage{amsmath,amssymb}
\usepackage{amsthm}
\usepackage[]{harpoon}
\usepackage[active]{srcltx} 
\usepackage[]{amsfonts}
\usepackage[]{amssymb}
\usepackage[]{harpoon}
\usepackage{multirow}
\usepackage{epsfig}
\usepackage{caption}
\usepackage{epstopdf}
\usepackage{booktabs}
\usepackage{subfigure}
\usepackage{color}%

\usepackage{lineno,hyperref}
\usepackage{showkeys}

\def\bbf{\mathbf{f}}
\def\bg{\mathbf{g}}
\def\bh{\mathbf{h}}

\def\br{\mathbf{r}}

\def\bu{\mathbf{u}}
\def\bv{\mathbf{v}}
\def\bw{\mathbf{w}}

\def\bA{\mathbf{A}}
\def\bB{\mathbf{B}}

\def\bI{\mathbf{I}}

\def\bL{\mathbf{L}}

\def\n{\noindent}
\def\pt{\partial}

\def\Re{\mathbb{R}}

\def\f#1#2{\frac {#1}{#2}}
\def\f32{\frac 32}

\def\bga{\begin{array}}
\def\eda{\end{array}}

\def\De{\Delta}

\def\ev{\equiv}

\def\gm{\gamma}

\def\la{\lambda}

\def\al{\alpha}

\def\rw{\rightarrow}
\def\iy{\infty}
\def\Om{\Omega}

\def\dfr#1#2{\displaystyle{\frac{#1}{#2}}}

\begin{document}
\newtheorem{thm}{Theorem}[section]
\newtheorem{cor}[thm]{Corollary}
\newtheorem{lem}[thm]{Lemma}
\newtheorem{prop}[thm]{Proposition}
\theoremstyle{definition}
\newtheorem{defn}[thm]{Definition}
\newtheorem{asser}[thm]{Assertion}
\theoremstyle{remark}
\newtheorem{rem}[thm]{Remark}
\numberwithin{equation}{section}

\newcommand{\norm}[1]{\left\Vert#1\right\Vert}
\newcommand{\abs}[1]{\left\vert#1\right\vert}
\newcommand{\set}[1]{\left\{#1\right\}}
\newcommand{\Real}{\mathbb R}
\newcommand{\eps}{\varepsilon}
\newcommand{\To}{\longrightarrow}
\newcommand{\BX}{\mathbf{B}(X)}
\newcommand{\A}{\mathcal{A}}
\newcommand{\ts}{\textstyle}
\newcommand{\tg}{\mbox{\rm tg}}
\newcommand{\ctg}{\mbox{\rm ctg}}
\newcommand{\atg}{\tg^{-1}}
\newcommand{\actg}{\ctg^{-1}}
\newcommand{\asin}{\sin^{-1}}
\newcommand{\acos}{\cos^{-1}}
\newcommand{\dps}{\displaystyle}
\newcommand{\fnz}{\footnotesize}
\newcommand{\D}{\displaystyle}
\newcommand{\DF}[2]{\frac{\D#1}{\D#2}}
\newcommand{\scp}{\scriptstyle}

\title{One-sided GRP Solver and Numerical Boundary Conditions for compressible fluid flows}

\thanks{   }

\author{  Jiequan Li}
\address{Institute of Applied Physics and Computational Mathematics, Beijing;  Center for Applied Physics and Technology, Peking University, China; and State Key Laboratory for Turbulence Research and Complex System, Peking University, China}

\author{Qinglong Zhang} 
\address{Department of Mathematics and Statistics, Ningbo University,  China }

\email{
  Jiequan Li:  li\_jiequan@iapcm.ac.cn}
  \email{Qinglong Zhang: zhangqinglong@nbu.edu.cn}


\maketitle

\begin{abstract}
 In the computation of compressible  fluid flows, numerical boundary conditions are always necessary for all physical variables at computational boundaries while just  partial physical variables are often prescribed as physical boundary conditions.  Certain extrapolation technique or ghost cells are  often  employed traditionally for this issue but  spurious wave reflections often arise to cause numerical instability.  In this paper, we associate this issue with the one-sided generalized Riemann problem (GRP) solver motivated by the accelerated piston problem in gas dynamics so that the extrapolation technique can be actually avoided. In fact, the compatibility arguments naturally requires to formulate the one-sided generalized Riemann problem and incorporate it into the numerical procedure of boundary conditions.  
 As far as the interaction of nonlinear waves with physical boundaries, such a one-sided GRP solver shows significant effects, as numerical experiments demonstrate, on avoiding  spurious wave reflections at the computational boundaries.
  \end{abstract}



\

\pagestyle{myheadings} \thispagestyle{plain} \markboth{Jiequan Li and Qinglong Zhang} {One-sided GRP and Numerical  Boundary Conditions}

\section{Introduction}\label{sec-Intro}

The issue on  boundary conditions for hyperbolic problems and particularly for compressible fluid flows is a classic topic and so is the corresponding numerical treatment.  There are a number of contributions via various approaches in literature, which   roughly  consist of  three types of concerns: physical considerations\cite{P-L-1992}, mathematical justifications (well-posedness arguments) \cite{Kreiss-book} and numerical treatment.   Physical considerations prescribe boundary data for a part of  physical variables based on specific problems such as the solid-wall boundary condition; mathematical 
 well-posedness justifies the validity of modelings subject to the prescribed boundary conditions; while numerical boundary conditions are prescribed for all physical variables so that discrete (approximate) equations can be implemented practically.  These concerns, though with different objectives, have the common goal on correctly describing the underlying problems, for which the  compatibility among the governing equations, prescribed boundary conditions and the initial data is a fundamental issue. 
  As far as the numerical treatment is concerned, extrapolation technique is often employed, particularly for high order accurate numerical methods. For example, 
in \cite{KP1,KP2} a lagrangian  interpolation is performed to achieve a second order accurate approximation to the boundary data in space.  However,  it just gives the first order accurate approximation in time.  Other works can be found, e.g. in \cite{Berger, FJ},  in the finite volume framework, and even in complex geometries \cite{KB}.   In \cite{P-L-1992}, characteristic method is used but restricted to first order accuracy for smooth flows.  Although the resulting schemes may be well-implementable, the corresponding validation is not clear both from rigorous mathematical analysis and numerical performance.  Improper extrapolation may lead to numerical instability such as spurious oscillations \cite{P-L-1992}.  From the viewpoint of numerical analysis, it is questionable whether the numerical boundary conditions are compatible with the discretized governing equations even though the underlying PDE models are well-posed. Hence it is worth addressing issue even though there are lots of studies available \cite{Berger,DDJ,DL1,DL3,FAH,FJ,Kreiss-book,KP1,KP2,KB,SP,TS1,TWSN}. 

We  associate this issue with the so-called {\em one-sided generalized Riemann problem (GRP)}. As motivation, we take a look at the initial-boundary value problem for the Burgers equation \cite{Bardos-1979,Bank-Ben-Artzi},
\begin{equation}
\begin{array}{ll}
u_t +(u^2/2)_x=0, \ \ \ &x\in (0,L), t>0,\\
u(x,0)=u_0(x), \ & x\in (0,L), 
\end{array} 
\label{eq:burgers}
\end{equation} 
and focus on the left boundary $x=0$. 
We assume that $u_0(x)\ev 1$ for example and inspect various situations upon the boundary requirement on $x=0$. There are three typical cases:
\begin{enumerate}
\item[(i)] $u(0,t)=a$, $0\leq a\leq 1$.  For this case, the solution contains a rarefaction wave
\begin{equation}
u(x,t) =\left\{
\begin{array}{ll}
a, & 0\leq x/t<a,\\
x/t, & a\leq x/t\leq 1,\\
1, &x/t>1. 
\end{array}
\right. 
\end{equation} 

\item[(ii)] $u(0,t)=b$, $b>1$. For this case, we have a shock solution, 
\begin{equation}
u(x,t) =\left\{
\begin{array}{ll}
b, & 0\leq x/t<(b+1)/2,\\
1, & x/t>(b+1)/2,
\end{array}
\right. 
\end{equation} 
 
\item[(iii)] $u(0,t)=c$, $c<0$.  There exists no physically admissible solution for such a case and so the boundary condition is not well prescribed. 
\end{enumerate}

This example shows the subtlety of nonlinear problems as investigated in \cite{Bardos-1979},  unlike linear hyperbolic problems. 
 In fact, for linear problems, physical boundary conditions depend on characteristic propagations.  While for nonlinear compressible fluid flows, many physical boundary conditions are prescribed upon surroundings and  cannot be even given {\em a priori}, such as the interaction of shock with solid  boundaries \cite{Ben-Dor-2007} and the solid body floating in the air \cite{DL3}. Numerically, situations become more complicated. First, numerical boundary conditions should be given for all variables in order to be suitable for the computation so that proper extrapolations have to be used.  For the strong interaction of waves with physical boundaries, the nonlinearity actually prevents the validity of extrapolations that may result in factitious phenomena. Second, high order approximations of boundary conditions are often made independently of  the discretization of the governing equations, which may lead to incompatibility and loss of accuracy.

The GRP formulated here is different from the traditional  GRP \cite{BenF, BenL1} and more suitably called {\em one-sided GRP}. It is an initial-boundary value problem rather than a purely initial value problem.  Correspondingly, a numerical method to solve this problem is called a {\em one-sided GRP solver.} Such a study  has two-fold goals:  (i) It  is used to the compatibility of prescribed boundary condition with the governing equations; (ii)  It proposes a family of high order numerical boundary conditions effectively compatible with the discretized governing equations.  In fact, accelerated piston problems \cite{CF, DL3} are the one-sided GRP  formulated here and 
they have been refined to put into the simulation of fluid flows with moving boundaries \cite{FAH, DL3}.   In \cite{DL1}, the similar idea was employed if the flow is smooth and consistent with the inverse Lax-Wendroff method \cite{TS1}.  Note that  since the two-stage fourth order framework \cite{Li1, LD, DL1} can be used to develop high order methods, there is no need to compute derivatives of order more than second.  Hence this paper just focuses on second order GRP solvers, which provides a reliable tool no matter whether the solution is smooth or not.

We organize this paper as follows. In Section \ref{sec:OS-GRP},  we formulate the one-sided GRP problem. In Section \ref{sec:OS-GD} we discuss numerical boundary conditions for compressible Euler equations via the passage of one-sided GRP.  We implement the resulting scheme in Section \ref{sec:num-0} and particularly display numerical results   in Subsection \ref{sec:Numerical} to demonstrate the performance.

 \section{High Order Numerical boundary conditions and one-sided GRP solver  for hyperbolic balance laws}\label{sec:OS-GRP} 
 
 In this section we formulate one-sided Riemann problems for hyperbolic balance laws in a general framework and discuss related numerical solvers for  the construction of numerical boundary conditions.  This is different from the classical initial-boundary value problem for compressible fluid flows in \cite{Li-Yu} since  initial and boundary values are  generally not compatible in a continuous way. Such a setting is  proposed for practical request. For example,  the flow is discontinuous at the reflection point of a shock on a solid boundary. Correspondingly, the associated initial-boundary value problem is formulated below as the one-sided generalized Riemann problem (OS-GRP). 
 
 Consider  hyperbolic balance   laws 
 \begin{equation}
 \begin{array}{ll}
 \bu_t+\bbf(\bu)_x=\bh(x,\bu), \ \ \ \  & x\in(0, L), t>0,\\
 \bu(x,0)=\bu_0(x), &  x\in(0, L),
 \end{array}
 \label{eq:ini-bdy}
 \end{equation}
 where $\bh(x,u)$ is a source term representing external forces or geometrical effects, $\bbf(\bu)$ is the flux function.  
 This system includes the compressible Euler equations we specified in the next section and many other models \cite{Dafermos-book}.
It is assumed to be hyperbolic in the sense that 
 the Jacobian $\bA(\bu)$ of $\bbf(\bu)$ has $m$ real eigenvalues $\la_k$ with a complete set of  associated eigenvectors $\br_k$, 
 \begin{equation}
 \bA(\bu)\br_k=\la_k \br_k,\ \ \ \ \la_1\leq  \cdots\leq \la_m. 
 \end{equation} 
 Each $\la_k$ is genuinely nonlinear or linearly degenerate in the sense of Lax \cite{Lax-1957}.
 
 We focus  on the  left boundary $x=0$. The right boundary $x=L$  is treated similarly.  We emphasize that the free boundary problem can be studied too \cite{DL3}. 
 On the boundary $x=0$, the data is imposed as 
 \begin{equation}
 \mathscr{B}\bu=\bg(t)\in \Sigma\subset \Re^{m-k(\bu)},
 \label{eq:bdc}
 \end{equation} 
 where the operator $\mathscr{B}: \Om\subset \Re^m \rightarrow \Sigma\subset \Re^{m-k(\bu)}$ projects the solution onto the boundary $x=0$, $0\leq k(\bu)\leq m$ is the number of negative eigenvalues.  If \eqref{eq:ini-bdy} is a linear problem, i.e., $\bbf(\bu)=\bA \bu$, $\bA$ is a constant matrix,  then the operator  $\mathscr{B}$ can be expressed in the matrix form 
 \begin{equation}
\bB \bu =\bg(t),  \ \ \ \  rank (\bB)=m-k,
\label{eq:bdc-linear} 
 \end{equation} 
 if $\bA$ has $m-k$ positive eigenvalues, where $\bB$ is a $\ell\times m$ matrix.
 For nonlinear problems the eigenvalues depend on the solution $\bu$ and thus   the integer $k(\bu)$ may vary depending on the solution $\bu$ too. Hence the precise meaning of the operator is determined together with the solution of \eqref{eq:ini-bdy}, as shown for the Burgers equation.

 Denote $\gm: \Re^m\rw \Re^m$ the trace operator on the boundary $x=0$,
  \begin{equation} 
\bu(x,t)|_{x=0}=\gm \bu(x,t).
 \end{equation} 
Then we propose the following assumption. 
\vspace{0.2cm} 

\n {\bf Assumption.}  The problem \eqref{eq:ini-bdy} --\eqref{eq:bdc} is well-posed at least locally so that 
\begin{equation}
\mathscr{B}(\gm \bu)=g(t)
\end{equation}
in some ``appropriate"  sense. 

\vspace{0.2cm}

This assumption is very ``rough" and understood in certain intuitive way.  A primary judgement of the well-posedness  boils down to the solvability of the following one-sided Riemann problem, while the dynamics is dependent on the one-sided generalized Riemann problem (GRP).  Numerically each component of $\bu$ should be given a value on the boundary so that the corresponding  numerical code can be implemented. Even with extrapolation, the approximation should be consistent with the solution of this one-sided GRP up to some desired accuracy order. 

In this section, we denote by $\bu_b(t)=\bu(0,t)$ the boundary value for the solution $\bu$, and $(\pt \bu/\pt t)_b(t)=\pt \bu/\pt t(0,t)$ the derivative of $\bu$ along the boundary $x=0$. 

\subsection{Linear equations with constant coefficients} \label{subsec:linear} 

Let's first get motivation from linear equations.  Consider with the assumption as above
\begin{equation} 
\bu_t +\bA\bu_x=\bh(\bu,t), \ \ \ x>0, t>0.\\
\label{eq:linear}
\end{equation} 
The characteristic decomposition tells that 
\begin{equation}
\dfr{\pt v_i}{\pt t}  +\la_i \dfr{\pt v_i}{\pt x}  =\bL_i \bh(\bu,x), \ \ \ \ \ i=1,\cdots, m, 
\label{eq:ch-linear}
\end{equation} 
where $v_i =\bL_i \bu$, $\bL_i$ is a left-eigenvector associated with the eigenvalue $\la_i$.  The solution formula is 
\begin{equation}
\begin{array}{ll}
v_i(x,t) &= v_i(x-\la_i t,0) \\
&\displaystyle +\int_0^t \bL_i\bh(\bu(x-\la_i(t-s),s),x-\la_i(t-s))ds=:K_i(x,t), \ \ \ \ i=0,\cdots,k.
\end{array}
\end{equation} 
 To obtain the solution values on the boundary $x=0$, we have from 
\eqref{eq:ch-linear}
\begin{equation}
v_i(0,t) = v_i(-\la_i t,0) +\int_0^t \bL_i\bh(\bu(-\la_i(t-s),s),-\la_i(t-s))ds=:K_i(0,t),  \ \ i=0,\cdots,k.
\end{equation} 
Hence the boundary value of $\bu$ can be  obtained by solving the following system
\begin{equation}
\begin{array}{l}
v_i= \bL_i \bu =K_i(t),  \ \ \ i=1,\cdots, k,\\
\bB\bu(0,t)=\bg(t).
\end{array}
\label{eq:linear-bd}
\end{equation} 
Indeed, the well-posedness of \eqref{eq:linear} depends on the solvability of \eqref{eq:linear-bd}. That is, 
\begin{equation}
\mbox{rank}\{\bL_1,\cdots, \bL_k, \bB_1, \cdots, \bB_\ell\}=m,
\end{equation} 
where $\bB_i$, $i=1,\cdots, \ell$, are  the row vectors of the matrix $\bB$. 
Such a solution formula in turn helps to develop high order schemes. We can refer to  \cite{DL1, TS1}  and next sections for the practical implementation in gas dynamics, corresponding to the acoustic case of one-sided GRP problem.  

\vspace{0.2cm} 

The above discussion is of course made  in the theoretical viewpoint. Numerically, we implement at each time level $t=t_n$, as follows. 
\begin{enumerate} 
\item[(i)] {\bf First order approximation.} The initial data is assumed to be constant  $\bu_R$.  Then we  derive all components of $\bu_b$  by solving the following system  
\begin{equation}
\begin{array}{l}
\bB u_b(t_n) =g(t_n)\\
v_i(0,t_n) =(v_i)_R, \ \ i=1,\cdots, k. 
\end{array}
\end{equation} 
Obviously, this is exactly the same as the usual extension from the neighboring  interior point using the characteristic method. 
\vspace{0.2cm}

\item[(ii)] {\bf Second order approximation.}  As high order approximations are concerned, we not only need to know the value in the first order approximation, but we have to approximate the value $(\pt \bu/\pt t)_b(t_n)=(\pt \bu/\pt t)(0,t_n)$ as well.  We denote by $\bu_R':=\bu_0' (0+0)$ and subsequently $\bv_R'=\bv'_0(0+0)$. Then we have
\begin{equation}
\begin{array}{l}
\bB  (\pt \bu/\pt t)_b(t_n) =g'(t_n)\\
(\pt v_i/\pt t)_b  (t_n) =-\la_i (v_i)_R'+\bL_i \bh(\bu_0(0),0), \ \ \ \ i=1,\cdots, k. 
\end{array}
\end{equation} 
Solving this system yields the value $(\pt \bu/\pt t)_b(t_n)$. This second order approximation shows clearly that the source term $\bh$ is input into the numerical boundary condition, unlike some direct extrapolation technique. Moreover, this characteristic method allows to deal with discontinuities at the origin. 

\end{enumerate} 

Indeed, these two approximations correspond to the one-sided Riemann problem and one-sided generalized Riemann problem, respectively.

 \subsection{One-sided Riemann problem} 
 
 The one-sided Riemann problem is motivated from the piston problem \cite{CF} and formulated in \cite{DL3}. Here we formulate this problem for 
 hyperbolic conservation laws
  \begin{equation}
 \begin{array}{ll}
 \bu_t+\bbf(\bu)_x=0, \ \ \ \  & x\in(0, \iy), t>0,\\
 \bu(x,0)=\bu_R, &  x\in(0, \iy).
 \end{array}
 \label{eq:RP}
 \end{equation}
The  boundary data is prescribed as 
 \begin{equation}
 \mathscr{B}(\bu)=\bv^*\in \Re^{m-k}, 
 \label{eq:b-data}
 \end{equation} 
for some $k\geq 0$, where the operator $\mathscr{B}$ prescribes certain physically meaningful values to partial state components. 
Corresponding to \eqref{eq:ini-bdy}-\eqref{eq:bdc}, $\bu_R =\bu_0(0+0)$ and $\bv^*=\bg(0)$. 

In order to solve this problem, we can mimick the method for the standard Riemann problem in the state space \cite{CH, Lax-1957, Toro-book}. At least for Euler equations, we will show how to solve it in the next section. 
The solvability of this one-sided Riemann problem depends on the compatibility of  the prescribed boundary conditions with the initial data. Generally speaking, as shown for the Burgers equation, this problem may not have to be well-posed. Hence this one-sided Riemann problem plays a  role in checking whether the boundary conditions are correctly prescribed. 

Another role of the one-sided Riemann problem is to supplement all state variables for the practical calculation because the boundary conditions just prescribe partial components of them.  For instance, consider the linear case, as indicated in \eqref{eq:bdc-linear}, with $m-k$ characteristics leaving the boundary $x=0$ so that the rank of the boundary operator $B$ is $m-k$.  Then we use the characteristic decomposition to obtain other $k$ equations, as shown above.

Assume that we are able to solve this problem  and obtain the solution $\bu(x,t)$  with the  trace on the boundary $x=0$ such that, 
\begin{equation}
\begin{array}{ll}
 \bu(x,t)|_{x=0} =\bu^*,\\
 \mathscr{B}(\bu^*) =\bv^*.
\end{array} 
\end{equation} 
Then it is  necessary  to check  whether there are exactly $m-k$ characteristics leaving from the boundary $x=0$, similar to the linear case. 
\begin{equation}
0<\la_{k+1}(\bu^*)\leq \cdots \leq \la_m(\bu^*). 
\end{equation} 
That is, the dimension of manifold $dim \{\mathscr{B} (\bu) =\bv^*\}=m-k.$ Just like the case for the Burgers equation, this is not necessary true. Hence  the solvability of one-sided Riemann problem is a necessary to judge the well-posedness of initial-boundary value problem for \eqref{eq:bdc-linear}.

 \subsection{One-sided GRP}
 As \eqref{eq:ini-bdy} includes a source term or/and the initial condition is not uniform (typically consists of piecewise polynomials), one has to consider a one-sided generalized Riemann problem (GRP). From the numerical point of view, one needs to have high order accurate prescription of all components of $\bu$ on the boundary $x=0$ as well as the construction of spatial variation near the boundary when high order methods are sought. 
 
 For completeness, the one-sided generalized Riemann problem (GRP) is reformulated here as the initial and boundary value problem, 
 \begin{equation}
 \begin{array}{ll}
 \bu_t +\bbf(\bu)_x = \bh(x,\bu), \ \ \ & x\in (0,\iy), t>0,\\
 \bu(x,0)=\bu_0(x),  & x\in (0,\iy),\\
 \mathscr{B}(\gm\bu)(0,t)=\bv^*(t)\in \Re^{m-k(\bu)}, \ \ \ &t>0,
 \end{array}
 \label{eq:os-grp} 
 \end{equation}
 where $\bu_0(x)$ is smooth, and $\bv^*(t)$ is measurable.  This is associated with the one-sided Riemann problem above. Similar to the interrelation between  the standard generalized Riemann problem and the associated Riemann problem, we have the following proposition \cite{BenL1}. 
 
 \begin{prop} Assume that \eqref{eq:os-grp} is well-posed and let $\bu(x,t)$ be its solution.  Denote that $\bu^A(x/t;\bu_R,\bv^*)$ be the solution of the associated one-sided Riemann problem \eqref{eq:RP}--\eqref{eq:b-data}. Then for every direction $\al=x/t>0$, 
 \begin{equation}
 \lim_{t\rw o^+} \bu(\al t,t) =\bu^A(\al; \bu_R,\bv^*).
 \end{equation}
 This implies the wave configuration of \eqref{eq:os-grp} is the same as that of \eqref{eq:RP}-\eqref{eq:b-data} asymptotically. 
 \end{prop}

Note that we assume that  the associated one-sided Riemann problem is uniquely solvable. 
Since the current paper is mainly concerned with a numerical algorithm for high order numerical boundary conditions, we leave aside  for the moment the investigation of the rigorous mathematical theory.

 \subsection{One-sided Riemann solver and one-sided GRP solver} 
 So-called solvers refer to the processes numerically solving the corresponding problems.  Standard numerical Riemann solvers can be found in \cite{Toro-book} and the generalized Riemann problem (GRP) solver in \cite{BenF, BenL1}.   The one-sided solvers proposed here are associated with the Riemann solver \cite{Toro-book}  and the GRP solver \cite{BenF,BenL1}.
These solvers aims  (i) to provide all physical variables on the boundary $x=0$; (ii)  to apply the inverse GRP to inspect the interaction of boundary and initial data.  
 
Note that the boundary value $\bu(0,t)$ that we obtain is not necessary to be continuous with the initial data $\bu(x,0)$ at the origin $(x,t)=(0,0)$. 
If so, the solution is discontinuous.  Such observation is heuristic when dealing with the interaction between shocks and solid boundaries. 
Besides, such a process provides several indications: 
\begin{enumerate}
\item[(i)] The compatibility of the resulting boundary data $\bu(0,t)$ and the initial data $\bu(x,0)$ determines the regularity of flows (solutions) around the origin locally.  The one-sided Riemann solution is a key clue to the well-posedness.  The one-sided Riemann solver aims to find $\bu_b(0)$ numerically. 

\item[(ii)]  The one-sided GRP solution depends on the associated Riemann solution, and  the corresponding GRP solver aims to find the value $(\pt\bu/\pt t)_b(0)$ and  helps to build high order numerical schemes.

\end{enumerate} 

%
%
%
%

\subsection{One-sided GRP solver in two dimensions} 
We extend the one-sided GRP solver to two dimensions in this part.  Suppose we have a boundary $\mathscr{L}: \Gamma(x,y)=0$ which is independent of time.  For the boundary conditions that depend on the time such as a piston problem, we refer to \cite{DL3} for the associated GRP solver. Our strategy includes the following steps: we first solve a normal one-sided Riemann problem at any fixed point on the boundary $\mathscr{L}$ along the normal direction, namely,
\begin{equation}
\begin{array}{ll}
\displaystyle \frac{{\rm} \partial {\bf u}}{{\rm} \partial t}+\frac{{\rm} \partial f({\bf u})}{{\rm} \partial x}+\frac{{\rm} \partial g({\bf u})}{{\rm} \partial y}=0,& \Gamma(x,y)>0,t>0, \\[9pt]
{\bf u}(x,y,0)={\bf u}_{R}(x,y),\ \ \ \ \ \quad & \Gamma(x,y)>0,\\[3pt]
 \mathscr{B}(\gamma {\bf u})(x,y, t)={\bf v}^{*}(x,y,t), &\Gamma(x,y) =0,  t>0,
\end{array}
 \label{eq:os-grp-2d} 
\end{equation}
where the boundary value is prescribed to be ${\bf v}^{*}(x,y,t)$. Denote by ${\bf n}(x,y,t)$ the unit normal vector of $\mathscr{L}$. For the presentation simplicity, the boundary is set  along the $y$-axis, thanks to the Galilean invariance for fluid dynamical systems. Then \eqref{eq:os-grp-2d} can be transformed  to solving the following normal generalized Riemann problem along the $y$-axis,
\begin{equation}
\begin{array}{ll}
\displaystyle \frac{{\rm} \partial {\bf u}^{gal}}{{\rm} \partial t}+\frac{{\rm} \partial f({\bf u}^{gal})}{{\rm} \partial x}+\frac{{\rm} \partial g({\bf u}^{gal})}{{\rm} \partial y}=0,& x>0,t>0, \\[9pt]
{\bf u}^{gal}(x,y,0)={\bf u}^{gal}_{R}(x,y),\ \ \ \ \ \quad & x>0,\\[3pt]
 \mathscr{B}(\gamma {\bf u}^{gal})(x,y, t)={\bf v}^{*,gal}(x,y,t), &x=0,  t>0,
\end{array}
 \label{eq:os-grp-2d-g} 
\end{equation}
where $\bu^{gal}$, $\bu_R^{gal}$ and $\bv^{*, gal}$ are  the Galilean transform of $\bu$, $\bu_R$  and ${\bf v}^{*}$, respectively.  After resolving $\bu^{gal}$, we transform back to obtain $\bu$.  The same as the 1-D case, we solve the normal conservation law at $(0,y^*)$
\begin{equation}
\begin{array}{ll}
\displaystyle \frac{{\rm} \partial {\bf u}^{N}}{{\rm} \partial t}+\frac{{\rm} \partial f({\bf u}^{N})}{{\rm} \partial x}=0, \\[9pt]
{\bf u}^{N}(x,t=0)={\bf u}_{R}^{gal}(0,y^{*}),\quad x>0,\\[3pt]
 \mathscr{B}(\gamma {\bf u}^{N})(0,t)={\bf v}^{*,gal}(0,y^{*}, 0),\quad  t>0,
\end{array}
 \label{eq:os-normal-grp} 
\end{equation}
to obtain the normal Riemann solution $\bu^N$.  Then  we solve the following IBVP,
\begin{equation}
\begin{array}{ll}
\displaystyle \frac{{\rm} \partial {\bf u}^{gal}}{{\rm} \partial t}+\frac{{\rm} \partial f({\bf u}^{gal})}{{\rm} \partial x}=-\left(\frac{{\rm} \partial g({\bf u})}{{\rm} \partial y}\right)^{N},\\[9pt]
{\bf u}^{gal}(x,y,0)={\bf u}_{R}^{gal}(x,y),\quad x>0,y\in \mathbb{ R},\\[3pt]
 \mathscr{B}(\gamma {\bf u}^{gal})(0,y, t)={\bf v}^{*,gal}(0,y, t),\quad y\in \mathbb{R}, t>0,
\end{array}
\end{equation}
to obtain 
\begin{equation}
\left(\dfr{\pt\bu^{gal}}{\pt t}\right)_b=\lim_{t\rw 0} \dfr{\pt\bu^{gal}}{\pt t}(0,y^*,t)
\end{equation} 
at $(0,y^*)$, 
where the term $\left(\frac{{\rm} \partial g({\bf u})}{{\rm} \partial y}\right)^{N} = \frac{{\rm \partial }g}{{\rm \partial} {\bf u}}({\bf u}^{N})\left(\frac{{\rm \partial}{\bf u}}{{\rm \partial }y}\right)^{N}$  is a fixed value with the instantaneous value ${\bf u}^{*}$ obtained from \eqref{eq:os-normal-grp} and $\left(\dfr{\partial \bf u}{ \partial y}\right)^{N}$ interpolated from $\bu_R^{gal}$, reflecting the tangential effect along the boundary \cite{LD}.  Then the 2-D one-sided  GRP solver follows exactly the same as the 2-D GRP solver, one can find more details in \cite{LD}. 


\subsection{High order numerical boundary conditions} Once the one-sided GRP solver is available, the boundary data can be approximated with second order accuracy and the boundary volume can be dealt with as the ordinary control volume. That is, if at moment $t=t_n$, $\bu(0,t_n)$ and $(\pt \bu/\pt t)(0,t_n)$ are known, then the boundary flux is approximated in a common way,
\begin{equation}
\begin{array}{l}
\bu(0,t_n+\De t/2) =\bu^N +\dfr{\De t}2 \left(\dfr{\pt\bu^{gal}}{\pt t}\right)_b,\\[3mm]
\dfr{1}{\De t} \int_{t_n}^{t_{n+1}} \bbf(\bu(0,t))dt =\bbf(\bu(0,t_n+\De t/2)) +\mathcal{O}(\De t^2).
\end{array}
\end{equation} 
 Furthermore the integral of source term can be evaluated using the interface method too, 
\begin{equation} 
\begin{aligned}
\dfr{1}{\De t\De x}  \int_{t_n}^{t_{n+1}}\int_0^{\De x}  \bh(\bu(x,t))dx dt = \dfr 12( \bh(\bu(0,t_n+\De t/2)) & + \bh(\bu(\De x,t_n+\De t/2)) ) \\
& + \mathcal{O}(\De t^2+\De x^2). 
\end{aligned}
\end{equation}
 Analogously, we deal with the 2-D case.

 \section{Application to Gas Dynamical Systems} \label{sec:OS-GD}
 In this section we discuss the one-sided Riemann problem (RP) and the one-sided generalized Riemann problem (GRP) for gas dynamical systems. We first discuss one-dimensional case in the form \eqref{eq:ini-bdy}  with 
\begin{equation}
\bu =(\rho, \rho v, \rho E)^\top,  \ \ \ \bbf(\bu) =(\rho v, \rho v^2+p, v(\rho E+p))^\top,
\label{eq:Euler-1}
\end{equation}
 and $\bh(\bu,x)$ is a problem-dependent source term.  In particular, for nozzle flows $\bh(\bu,x)$ takes the form
 \begin{equation}
 \bh(\bu,x)=-\frac{a'(x)}{a(x)}\left(
\rho v,
\rho v^2,
v(\rho E+p) \right)^\top.
\label{eq:euler-source}
 \end{equation} 
 where $\rho, v, p$ are the density, velocity and pressure of the fluids, respectively. $a(x)$ is the cross-section area of the duct.  $E=\frac{v^2}{2}+e$ is the total energy, the internal energy $e$ is given by the equation of state (EOS) $e=e(\rho,p)$.   Note that \eqref{eq:euler-source} also includes the case of  radially symmetric flows \cite{LLS}. The reactive Euler flows \cite{TWSN} can be treated similarly.

\vspace{0.2cm} 

System \eqref{eq:Euler-1} has three eigenvalues 
\begin{equation}
\la_-=v-c, \  \ \ \la_0=v, \ \ \ \la_+ =v+c,
\end{equation} 
where $c$ is the local sound speed. All other properties of this system can be found in any textbook about gas dynamics, e.g. \cite{CF, CH, BenF, Toro-book}. 

{\bf The one-sided Riemann solver for Euler equations.}  As we pointed out in the last section, the one-sided Riemann problem plays  a role in the justification of local well-posedness besides its numerical value. 
   This problem is formulated as 
\begin{equation} 
\begin{array}{l}
\bu_t+ \bbf(\bu)_x=0, \ \ \ x>0,t>0,\\
\bu(x,0)=\bu_R, \ \ x>0,\\
\mathscr{B}\bu(0,t)=\bw_b. 
\end{array}
\label{eq:osrp}
\end{equation} 
The method solving this problem \eqref{eq:osrp}  follows the one for the classical Riemann problem. We fix the wave curve $W_R$ associated with $\la_+=v+c$ from the state $\bu_R$ in the phase space, $(\rho, v, p)$--space, and then investigate the solvability for the prescribed data $\bw_b$.  It is easily checked that the solvability of such a problem is up to the following two conditions:

\begin{enumerate}
\item[(i)] There is an intersection  point $\bu^*$  of  $W_R$ and $\mathscr{ B}\bu= \bw_b$; 

\item[(ii)] The dimension $dim\{\mathscr{B}\bu=\bw_b\} = \#\{\la_i(\bu^*)>0, i=-,0,+\}$.   

\end{enumerate} 

\begin{figure}[h]
\subfigure{
\begin{minipage}[t]{0.425\textwidth}
\centering
\includegraphics[width=\textwidth]{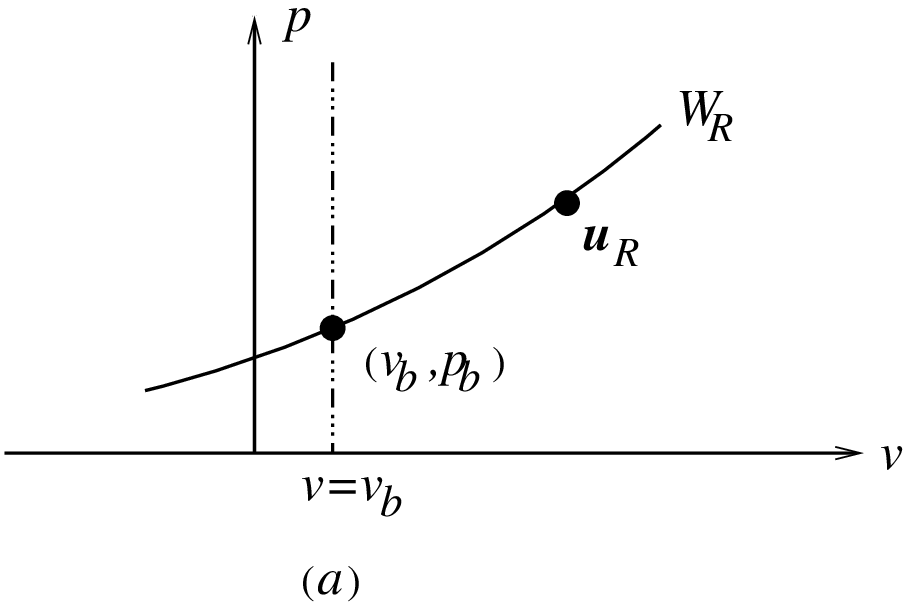}
\end{minipage}
}
\subfigure{
\begin{minipage}[t]{0.425\textwidth}
\centering
\includegraphics[width=\textwidth]{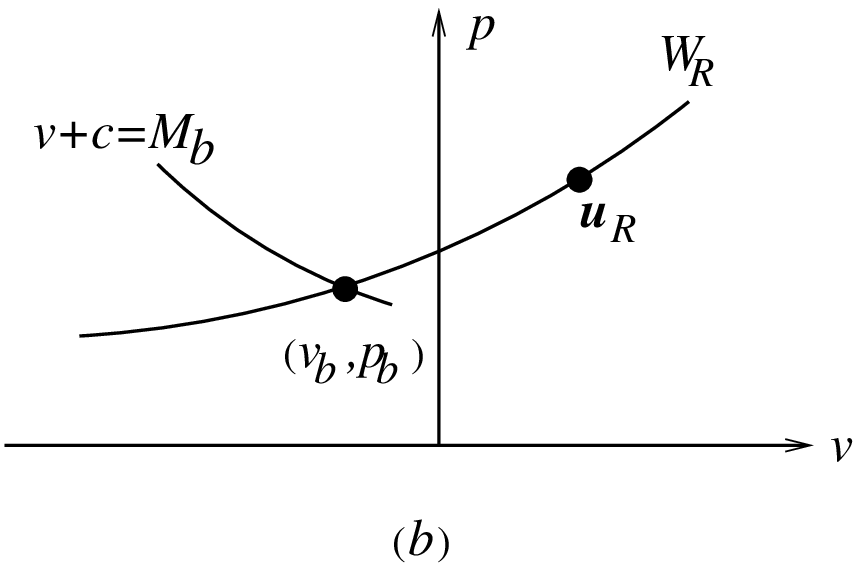}
\end{minipage}
}
\caption{ \scriptsize The one-sided Riemann problem for two typical boundary conditions: (a) Prescribed velocity $v_b$. (b) Given upstream Mach number $M_b$.}
\label{fig.1.}
\end{figure}

The following are two typical examples, see Fig. 3.1. 
\begin{enumerate}
\item[(i)] {\bf Prescribed velocity $v_b$. } Then we  find a point on $W_R$ so that $p_b$ is fixed. Certainly as $v_b>0$, $\la_0=v_b>0$ and $ \la_+=v_b+c_b>0$ so that one additional condition is needed.  For instance, we can supplement the gas density $\rho_b$ as the  gas property on the boundary $x=0$. However, as $v_b<0$, we need to check the Mach number $M_b=|v_b|/c_b$. If $M_b>1$, the boundary condition is not suitably prescribed. 

\vspace{0.2cm} 

\item[(ii)] {\bf Given  upstream Mach number $M_b$.  } The given value  $M_b=v_b/c_b$ actually implies 
\begin{equation}
v+c=M_b. 
\end{equation} 
We look for its intersection point with $W_R$ to find $(v_b,p_b)$ and then $\rho_b$ using the equation of state (EOS). 
\end{enumerate} 

In summary, we can  investigate the one-sided Riemann problem to identify that whether the upstream flow is supersonic or not as well as  make  clear the correct prescription of boundary conditions.

\vspace{0.2cm}

{\bf One-sided GRP solver.}  The one-sided GRP solver serves to solve \eqref{eq:os-grp} numerically.  Assume that $\bu_0(x)$ is (or approximated by)  a smooth function with regular limiting values 
\begin{equation} 
\bu_R =\lim_{x\rw 0+0} \bu_0(x),  \ \ \ \ \bu_R' =\lim_{x\rw 0+0} \bu_0'(x). 
\end{equation}
Based on the corresponding one-sided Riemann problem, we can obtain the limiting value  $(\pt\bu/\pt t)_b(0)$ on the boundary $x=0$. 
Essentially there are two versions in analogy with the standard GRP solver:  {\em An acoustic version} and {\em  a nonlinear version.}

\begin{enumerate}
\item[(i)] {\em Acoustic GRP.}  As $\|\bu_b-\bu_R\|\ll 1$, we can use the acoustic approximation, i.e., the linear method in Subsection \ref{subsec:linear}. 

\item[(ii)] {\em Nonlinear GRP.}  As strong waves emit from the corner $(0,0)$ (i.e., $\|\bu_b-\bu_R\|\gg 1$), we have to develop a genuinely nonlinear GRP solver, similar to the standard GRP solver for general hyperbolic balance laws \cite{BenL1}

\end{enumerate} 

Specified to the Euler equations, the one-sided GRP solver is implemented as follows. 
\begin{enumerate} 
\item[(i)] Judge from the associated one-sided Riemann solution whether there emit strong waves in order to determine to use the acoustic or nonlinear GRP solver. 

\item[(ii)] The acoustic GRP solver is the same as the linear case above. 

\item[(iii)] The nonlinear GRP solver consists of two cases: {\em  a supersonic upstream flow} and {\em a subsonic upstream flow.}\
\begin{enumerate}
\item {\em A supersonic upstream flow.}  All conditions are given at boundary $x=0$. 

\item {\em A subsonic upstream flow.} We apply the same procedure of the standard GRP solver \cite{BenLW} and naturally derive the one-sided relation
\begin{equation}\label{3.7}
\displaystyle a_{R}\left(\frac{\pt  v}{\pt t}\right)_b+b_{ R}\left(\frac{\pt p}{\pt t}\right)_b=d_{ R},
\end{equation}
where   the coefficients $a_{ R}, b_{ R}$ and $c_{ R}$ are fully determined by the values ${\bf u}_{R}(0), {\bf u}_{R}^{*}$ and the slope value ${\bf u}_{R}^{'}(0)$, the detailed expressions can be found in \cite{BenL1}.

We are in position to compute the partial derivative values $(\partial v/\partial t)_b$ and $(\partial p/\partial t)_b$ from \eqref{3.7}. 
If the boundary condition is given as $v_b(t)=g(t)$ and subsequently  $(\partial v/\partial t)_b=g'(t)$, then $(\partial p/\partial t)_b$ follows by   the linear relation \eqref{3.7}. As for the density derivative $(\partial \rho/\partial t)_b$, we have
\begin{equation}\label{3.8}
\left(\frac{\pt \rho}{\partial t}\right)_b=\frac{1}{(c_R^{*})^{2}}\left(\frac{{ \pt}p}{\partial t}\right)_b
\end{equation}
on the boundary from the EOS. Here $c_R^{*}$ is the local sound speed.

If the upstream boundary condition is given in terms of Mach number $M_b(t)=g(t)$, then one has 
\begin{equation}\label{3.9}
\left(\frac{\pt  v}{\pt t}\right)_b + c_p \left(\frac{\pt  p}{\pt t}\right)_b+c_\rho \left(\frac{\pt  \rho}{\pt t}\right)_b=g'(t),
\end{equation} 
where $c_p=\frac{\pt c}{\pt p}$ and $c_\rho=\frac{\pt c}{\pt \rho}$.  This, together with \eqref{3.7} and the relation \eqref{3.8}, provides the boundary condition. 

\end{enumerate} 

\end{enumerate} 

{\bf 2-D one-sided GRP solver for Euler.}  The 2-D compressible Euler equations can be written as 
\begin{equation}
\begin{array}{ll}
\displaystyle \frac{{\rm} \partial {\bf u}}{{\rm} \partial t}+\frac{{\rm} \partial \bbf({\bf u})}{{\rm} \partial x}+\frac{{\rm} \partial \bg({\bf u})}{{\rm} \partial y}=0,\\[9pt]
{\bf u}=\left(
\begin{array}{ccc}
\rho \\
\rho v^{x} \\
\rho v^{y} \\
\rho E 
\end{array}\right),
\quad 
\bbf({\bf u})=\left(
\begin{array}{ccc}
\rho v^{x} \\
\rho (v^{x})^2+p \\
\rho v^{x}v^{y} \\
v^{x}(\rho E+p) 
\end{array}\right),
\quad 
\bg({\bf u})=\left(
\begin{array}{ccc}
\rho v^{y} \\
\rho v^{x}v^{y} \\
\rho (v^{y})^2+p\\
v^{y}(\rho E+p) 
\end{array}\right),
\end{array}
\label{eq:2d-euler}
\end{equation}
where $\rho, v^{x}, v^{y}, p$ and $E$ represent the density, $x-$velocity, $y-$velocity, pressure and total energy, respectively. The  2-D one-sided GRP solver is the practical combination of the above 1-D one-sided GRP solver and the 2-D GRP solver \cite{BenL1, DL1}.  A key point is that  the transversal effect is included in the solver  development \cite{LL}.

\section{Implementation of the one-sided GRP scheme} 
\label{sec:num-0}

\subsection{ Brief summary of the one-sided GRP scheme} 
So far, the one-sided GRP solver is developed to suit the GRP scheme near the boundary.  The boundary control volume is then treated the same as the interior control volumes in the finite volume framework. We discretize the domain by equally computation mesh size $\Delta x=x_{j+\frac12}-x_{j-\frac12}$ and set $I_j=(x_{j-\frac12},x_{j+\frac12})$.  The cell $I_{0}=(x_{-\frac12},x_{\frac12})$ represents the left boundary cell centered at $x_0$ and the cell $I_{M}=(x_{M-\frac12},x_{M+\frac12})$ represents the right boundary cell centered at $x_M$, as shown in Fig. 4.1. 

\begin{figure}[htbp]
\centering
\includegraphics[width=0.6\textwidth]{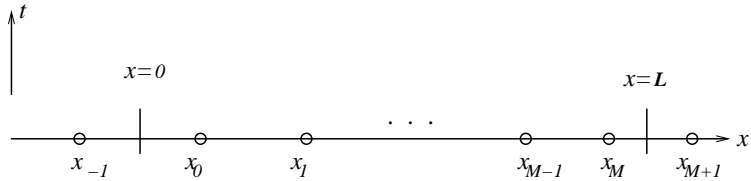}
\caption{\scriptsize The computational domain $(0,L)$. Set $x_0=\Delta x/2$ and $x_M=1-\Delta x/2$, The one-sided GRP is solved at the boundaries $x=0$ and $x=L$, respectively.}
\label{fig.2}
\end{figure}

Since  a standard finite volume method, such as the GRP method in \cite{BenLW}, can be applied over computational cells $I_j (j= 1,...,M-1)$ in the interior domain, we only focus on the boundary cell $I_0$. The one-sided GRP scheme at the boundary cell $I_0$ assumes the piecewise linear data 
\begin{equation}
{\bf u}(x,t_{n})={\bf u}_{0}^{n}+\sigma_0^{n}(x-x_0),\quad x\in(x_{-\frac12},x_{\frac12}).
\label{eq:initial}
\end{equation}
The vector $\sigma_{0}^{n}$ is the constant slope of ${\bf u}(x,t_{n})$ over cell $I_{0}$ at time $t_n=n\Delta t, n\in {\bf N}$ with $\Delta t$ the time step size. To obtain the second order accuracy, the mid-point value is used
\begin{equation}
{\bf u}_{-\frac12}^{n+\frac12}= {\bf u}(x_{-\frac12},(n+1/2)\Delta t)
\end{equation}
in the resolution of numerical flux and the source term discretization. We apply the 1-D one-sided GRP solver in the following steps.

\noindent
{\bf Step 1.} Given the piecewise linear initial data \eqref{eq:initial},  approximate the mid-point value ${\bf u}_{-\frac12}^{n+\frac12}$ as follows,
\begin{equation}
{\bf u}_{-\frac12}^{n+\frac12}={\bf u}_{-\frac12}^{n}+\frac{\Delta t}{2}\left(\frac{{\rm \partial}{\bf u}}{{\rm \partial}t}\right)_{-\frac12}^{n}.
\label{eq:mid-point}
\end{equation}
The computation of $({\rm \partial}{\bf u}/{\rm \partial}t)_{-1/2}^{n}$ is the main ingredient of the one-sided GRP scheme. The value ${\bf u}_{-\frac12}^{n}$ is the local solution at $(x_{-\frac12},t_n)$ to the following one-sided Riemann problem:
\begin{equation}
\left\{
\begin{array}{lll}
\displaystyle \frac{{\rm} \partial {\bf u}}{{\rm} \partial t}+\frac{{\rm} \partial \bbf({\bf u})}{{\rm} \partial x}=0,\\[9pt]
\mathscr{B} \bu =\bw_b, \quad x=x_{-\frac12}, \\[6pt]
{\bf u}_{R} :={\bf u}_{0}^{n}-(x_{0}-x_{-\frac12})\sigma_{0}^{n}, \quad x>x_{-\frac12},
\end{array}\right.
\end{equation}
which can be solved by an exact or approximate one-sided Riemann solver \cite{Toro-book}. Here we apply the one-sided GRP procedure \eqref{3.7}-\eqref{3.9} to obtain the instantaneous value $(\partial {\bf u}/\partial t)_{-1/2}^{n}$ on the boundary $x=x_{-\frac 12}$, and then approximate ${\bf u}_{-1/2}^{n+\frac12}$ using \eqref{eq:mid-point}.

\noindent
{\bf Step 2.} Evaluate the next time values ${\bf u}_{0}^{n+1}$ by using the following formula
\begin{equation}
{\bf u}_{0}^{n+1}={\bf u}_{0}^{n}-\frac{\Delta t}{\Delta x}\left(\bbf({\bf u}_{\frac12}^{n+\frac12})-\bbf({\bf u}_{-\frac12}^{n+\frac12})\right)+\frac{\Delta t}{2}\left(\bh(x_{\frac12},{\bf u}_{\frac12}^{n+\frac12})+\bh(x_{-\frac12},{\bf u}_{-\frac12}^{n+\frac12})\right),
\end{equation}
where the source term $\bh(x,{\bf u})$ is discretized with the mid-point rule in time and the trapezoidal rule in space.

\noindent
{\bf Step 3.} In order to suppress local oscillations as discontinuities are present near the boundary, we update the 
slope $\sigma_{0}^{n+1}$ by using the following monotonicity algorithm limiter
\begin{equation}
\sigma_{0}^{n+1}= {\rm minmod}\left(\displaystyle \frac{{\bf u}_{\frac12}^{n+1,-}-{\bf u}_{-\frac12}^{n+1}}{\Delta x}, \displaystyle \frac{{\bf u}_{1}^{n+1}-{\bf u}_{0}^{n+1}}{\Delta x}\right).
\label{eq:minmod}
\end{equation}
More details about the minmod function can be found in \cite{BenLW,Toro-book}.
\vspace{0.2cm} 
 
In two-dimensional computations, we take  rectangular meshes $\cup \bI_{j,k}$, $j=0,..., M$, $k=0,..., N$,   as an example for simplicity, here $\bI_{j,k}{\tiny =}(x_{j-1/2},x_{j+1/2})\times(y_{k-1/2},y_{k+1/2})$ centered at the grid point $(x_j,y_k)$. The finite volume formula is applied over all cells $\bI_{j,k}$, 
\begin{equation}
{\bf u}_{j,k}^{n+1}={\bf u}_{j,k}^{n}-\frac{\Delta t}{\Delta x}\left(f({\bf u}_{j+\frac12,k}^{n+\frac12})-f({\bf u}_{j-\frac12,k}^{n+\frac12})\right)-\frac{\Delta t}{\Delta y}\left(g({\bf u}_{j,k+\frac12}^{n+\frac12})-g({\bf u}_{j,k-\frac12}^{n+\frac12})\right).
\end{equation}
The initial data at time $t=t_{n}$ is expressed  as  bilinear functions
\begin{equation}
{\bf u}(x,y,t_{n})={\bf u}_{j,k}^{n}+(\sigma_x)_{j,k}^{n}(x-x_{j})+(\sigma_y)_{j,k}^{n}(y-y_{k}),\quad j=0,1,...,M, \quad k=0,1,...,N.
\label{eq:2d-initial}
\end{equation}
The values ${\bf u}_{j+\frac12,k}^{n+\frac12}$ and ${\bf u}_{j,k+\frac12}^{n+\frac12}$ can be analytically derived by the resolution of a local quasi 1-D GRP solver at each interface.  The one-sided GRP solver is applied on the boundary.  Then we can take the same procedure as that for 1-D case to implement the finite volume scheme.

\begin{figure}[htbp]
\centering
\includegraphics[width=0.6\textwidth]{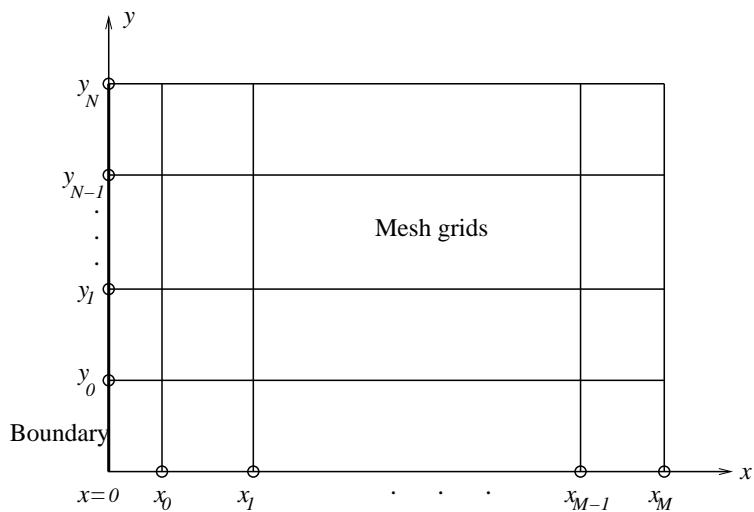}
\caption{\scriptsize The boundary is initially along the $y-$ axis. Set $x_0=\Delta x/2$ and the one-sided GRP is solved on the boundary $x=0$.}
\end{figure}

\subsection{Numerical Examples} \label{sec:Numerical}
We will present several numerical examples  to validate the performance as  the one-sided GRP solver is used. The examples include  the interaction of shocks with solid boundaries,  the radially symmetric flows, the nozzle flows, the Mach reflection of shock and the forward facing step problem.  

\begin{figure}[h]
\subfigure{
\begin{minipage}[t]{0.45\textwidth}
\centering
\includegraphics[width=1.0\textwidth]{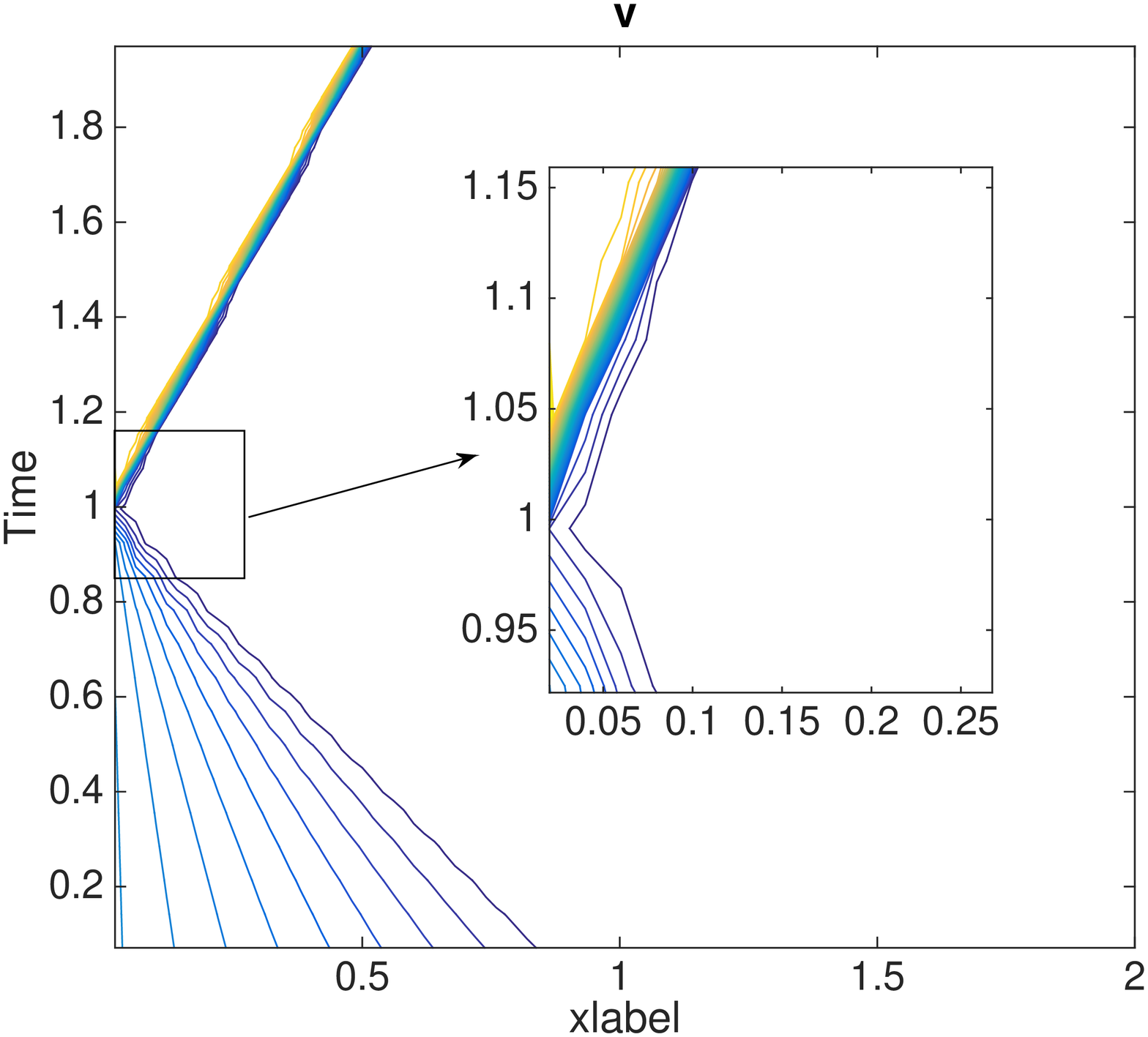}
\end{minipage}
}
\subfigure{
\begin{minipage}[t]{0.45\textwidth}
\centering
\includegraphics[width=1.0\textwidth]{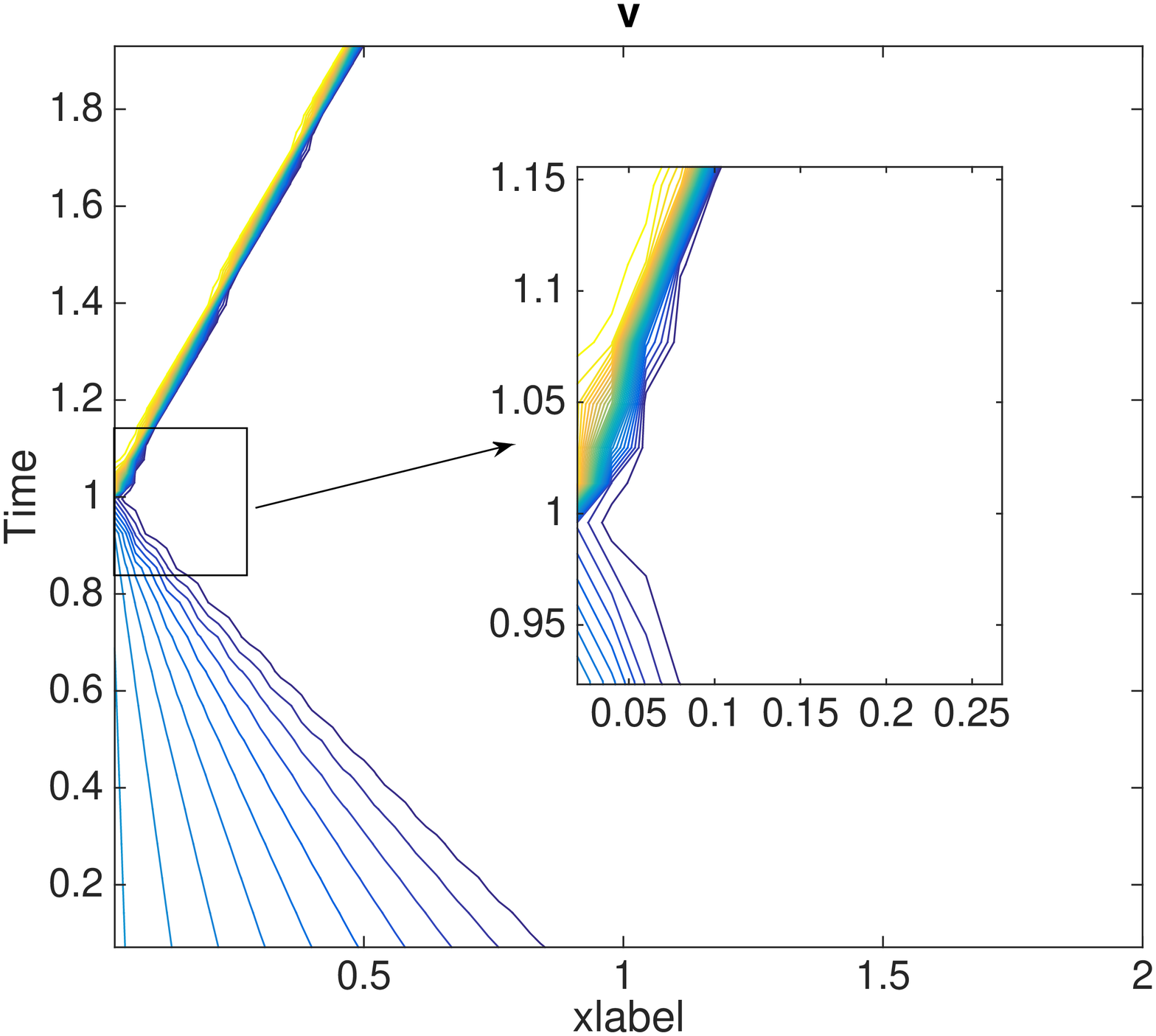}
\end{minipage}
}
\caption{\scriptsize The contours of  the solution $v(x,t)$ of the Burgers equation obtained by the one-sided GRP solver (left) and the traditional boundary condition treatment (right), respecticely. $100$ cells are used and  $30$  contours are drawn.}
\end{figure}
\vspace{0.2cm} 

\n {\bf Example 1. The scalar equation.}  We first use the Burgers equation to test the performance of the one-sided GRP solver. Consider the following initial-boundary value problem for the Burgers equation
\begin{equation}
\begin{array}{lll}
\displaystyle v_t+\left(\frac{v^2}{2}\right)_x=0, \ \quad x\in(0,2), t>0,\\[9pt]
v(x,0) = \left\{
\begin{array}{ll}
-x, \quad 0<x<1,\\
-1, \quad 1<x<2,
\end{array}
\right. \\[3mm]
v(0,t) = \left\{
\begin{array}{ll}
0,\quad 0<t<1,\\
2, \quad t>1.
\end{array}
\right.
\end{array}
\end{equation}
The solution $v(x,t)$ has an explicit formula
\begin{equation}
v(x,t) =\left\{
\begin{array}{lll}
\dfr{x}{t-1}, \ \ \ \ \ & 0<x<1-t, &0<t<1,\\
2, & 0<x<t/2,&t>1,\\
-1, & x>1-t, &0<t<1,\\
-1, & x>t/2,& t>1. 
\end{array}
\right. 
\end{equation} 
A compressible wave  propagates to the left and forms a shock at $(0,1)$ on the boundary.   As $t\geq 1$, a shock from $(0,1)$  propagates  to the right. We compute the solution using the GRP scheme with  the reflective boundary condition and the one-sided GRP solver, respectively.  The solution $v(x,t)$ is plotted from time $t=0$ to time $t=2$ in Fig. 4.3,  from which it is observed  that the one-sided GRP solver gives very sharp resolution of the singularity point $(0,1)$,  compared with the reflective boundary condition treatment.
\vspace{0.2cm} 

\n {\bf Example 2. A single shock interaction with a solid boundary}  We test the example that a single shock wave interacts with a solid wall to verify the numerical performance of the one-sided GRP solver. The computational domain is [0,10] where the boundary is at $x=0$. A left-propagating shock wave is initially positioned  at  $x=2$. We take $\gamma=1.4$ and the initial data is set to be
\begin{equation}
(\rho,v,p)(0,x)=\left\{
\begin{array}{ll}
(1.4,0.0,1.0),\quad 0\leq x\leq 2.0,\\[5pt]
(8.0,-8.25,116.5),\quad 2.0<x\leq 10.0.
\end{array}
\right.
\label{eq:initial-ex1}
\end{equation}

\begin{figure}[h]
\subfigure{
\centering
\includegraphics[width=0.45\textwidth]{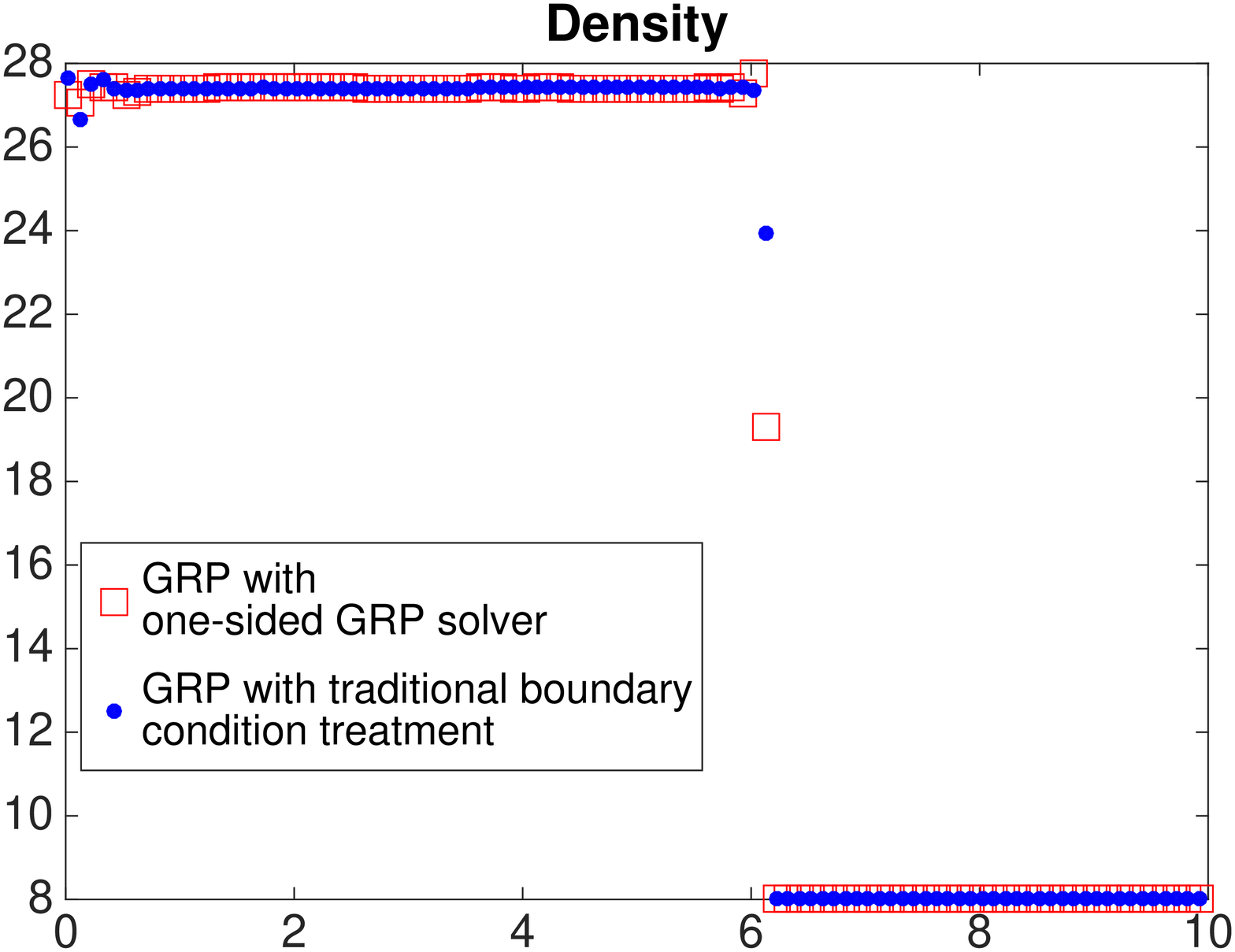}
}
\subfigure{
\centering
\includegraphics[width=0.45\textwidth]{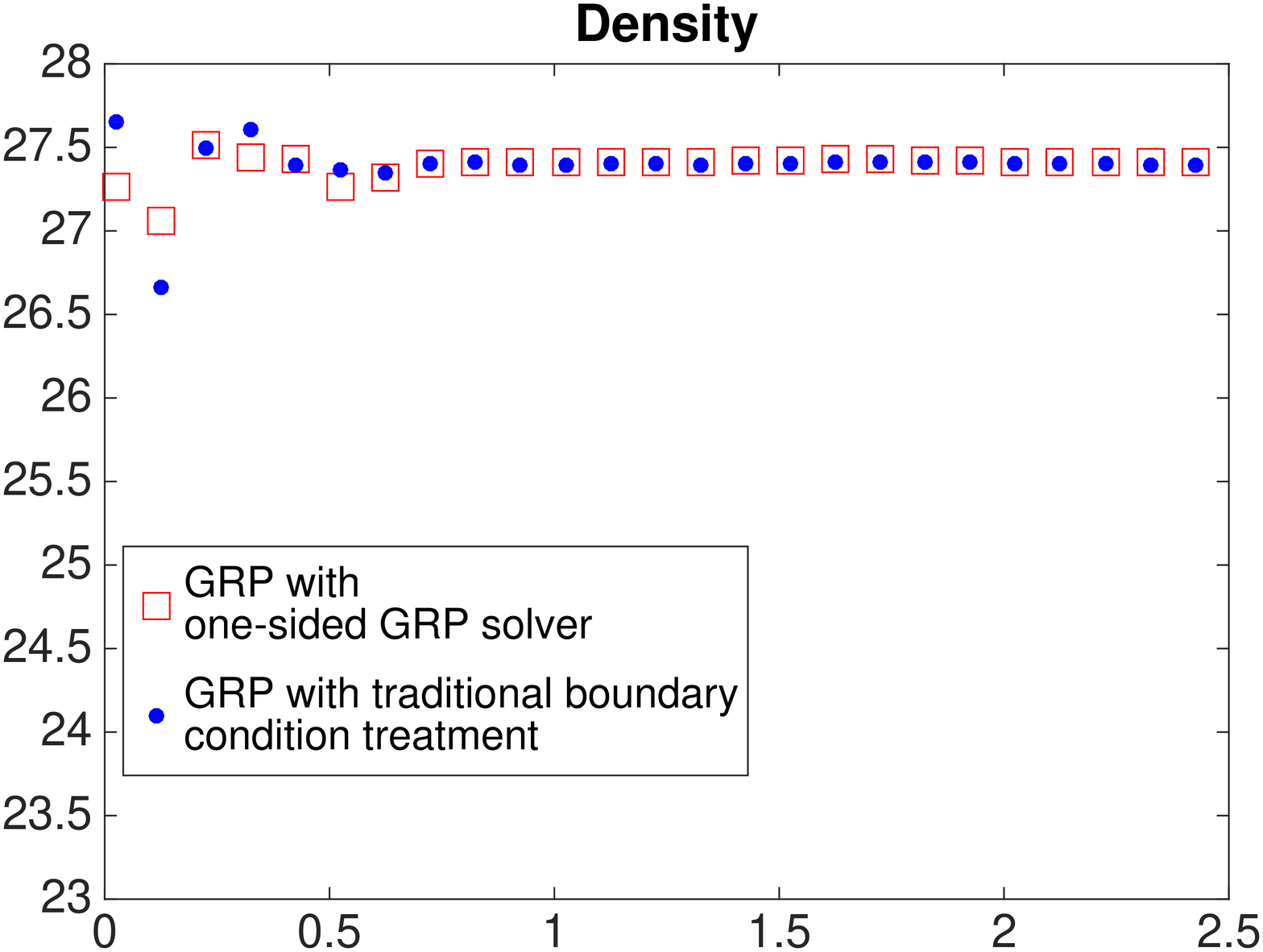}
}
\caption{\scriptsize A shock wave interacts with a solid wall. We compare the density profile obtained with the one-sided GRP solver (squares) with that obtained with the traditional reflective boundary condition (dots) with 400 cells (200 are shown).}
\end{figure}
\vspace{0.2cm} 

\begin{figure}[h]
\subfigure{
\begin{minipage}[t]{0.45\textwidth}
\centering
\includegraphics[width=1.0\textwidth]{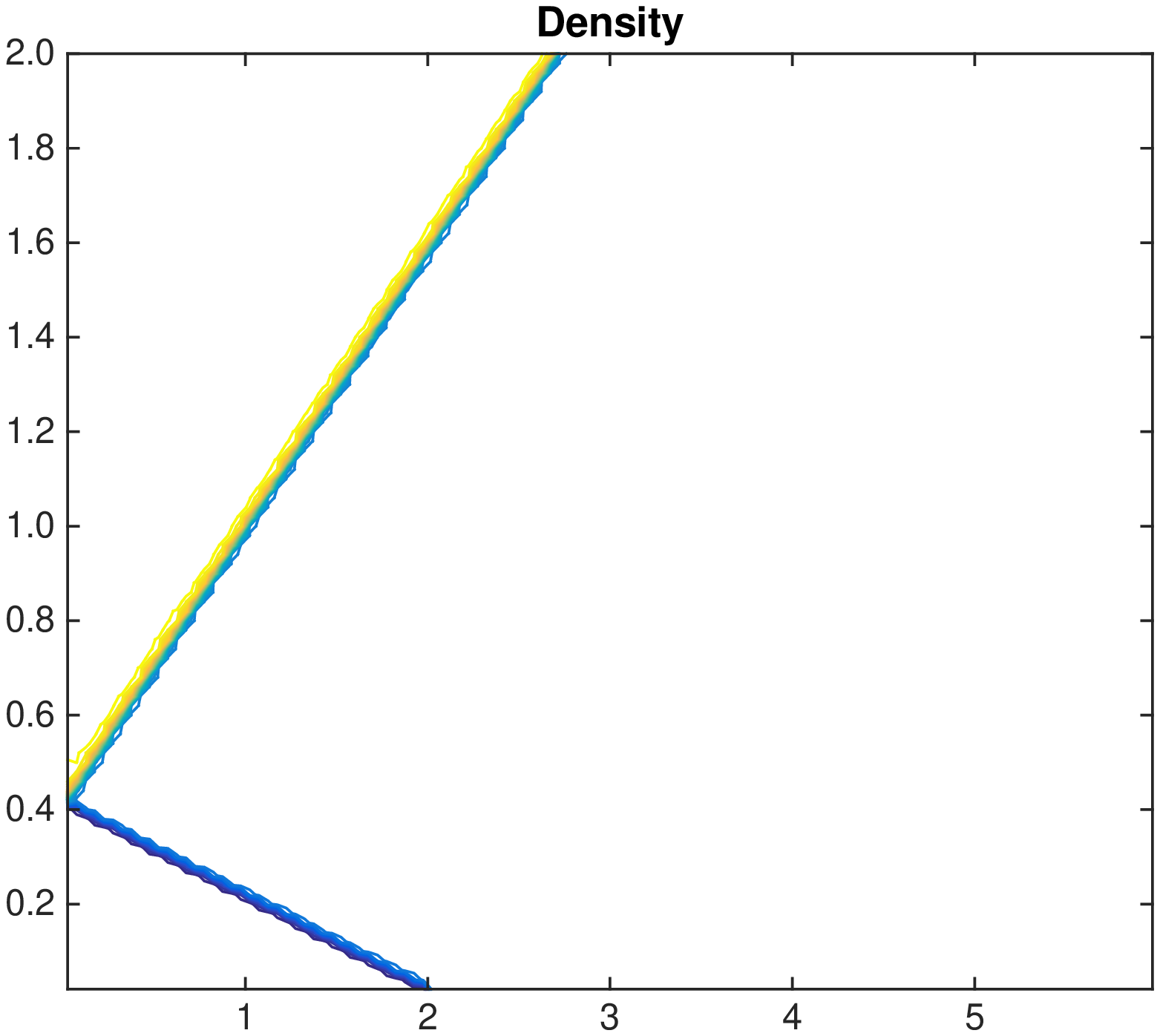}
\end{minipage}
}
\subfigure{
\begin{minipage}[t]{0.45\textwidth}
\centering
\includegraphics[width=1.0\textwidth]{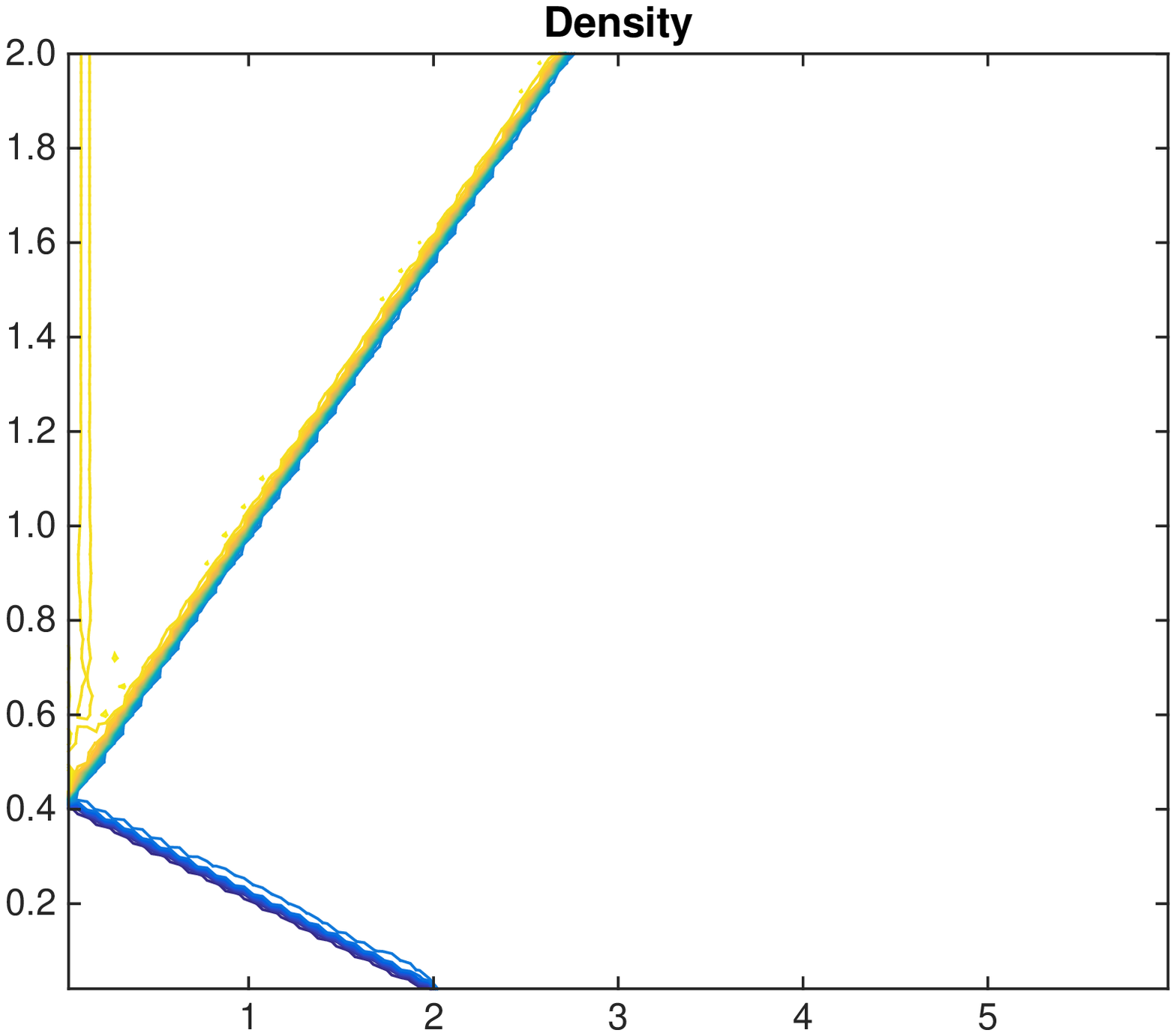}
\end{minipage}
}
\caption{\scriptsize The contours of the solution of Example $2$  obtained by the one-sided GRP solver (left) and the traditional boundary condition treatment (right). Thirty contours are drawn.}
\end{figure}
\vspace{0.2cm} 

A reflected shock wave is observed when the output time is set to $t=2.0$. We compare the results by two different boundary condition treatments: the traditional reflective boundary condition and the one-sided GRP solver. From Fig. 4.4, one can observe that the result obtained by the one-sided GRP solver is more stable and has less oscillations near the boundary. We further plot $30$ equally-distributed density contours from time $t=0$ to time $t=2$ at every  time interval $0.01$ in  Fig. 4.5, one can  see again  that the one-sided GRP solver gives very sharp resolution  at the interaction point of the shock with the boundary. 
\vspace{0.1cm} 

\noindent
{\bf Example 3. The Woodward-Colella problem.}  This is a classical interacting blast wave problem with the gas initially at rest and $\gamma=1.4$. The density is everywhere unit, the pressure is $p=1000$ for $0\leq x<0.1$ and $p=100$ for $0.9<x\leq 1.0$, while it is only $p=0.01$ for $0.1<x<0.9$. The solid-wall boundary conditions are prescribed at both ends. We compare the results of the reflective boundary condition treatment with that of the one-sided GRP solver. The CFL number is 0.6.
The output time is set to $t=0.038$. The numerical results for both boundary condition treatments are shown in Fig. 4.6 with 400 cells and 800 cells, respectively. It can be seen that the one-sided GRP solver is effective and robust for the blast wave problem.

\begin{figure}[h]
\subfigure{
\centering
\includegraphics[width=0.45\textwidth]{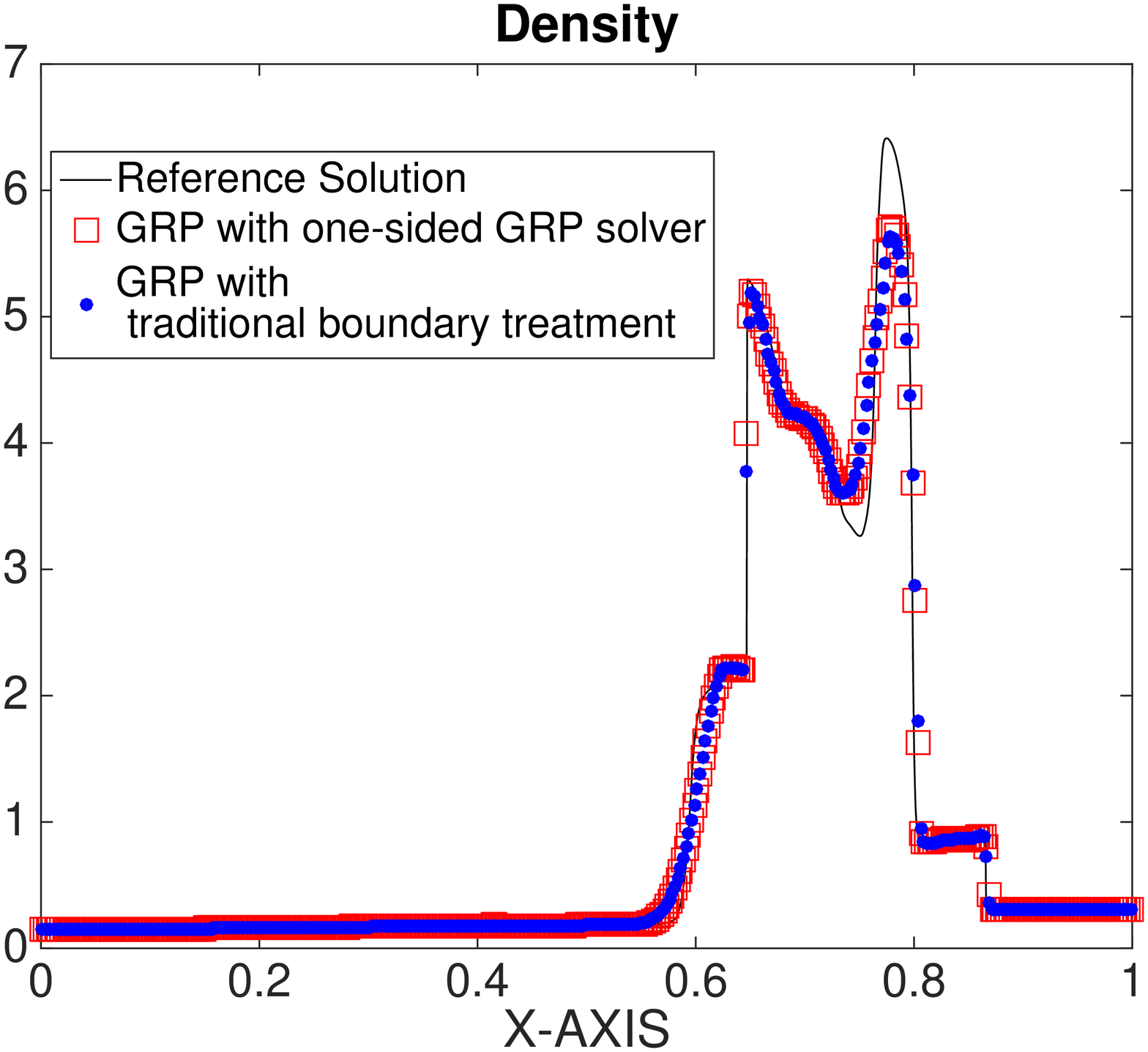}
}
\subfigure{
\centering
\includegraphics[width=0.45\textwidth]{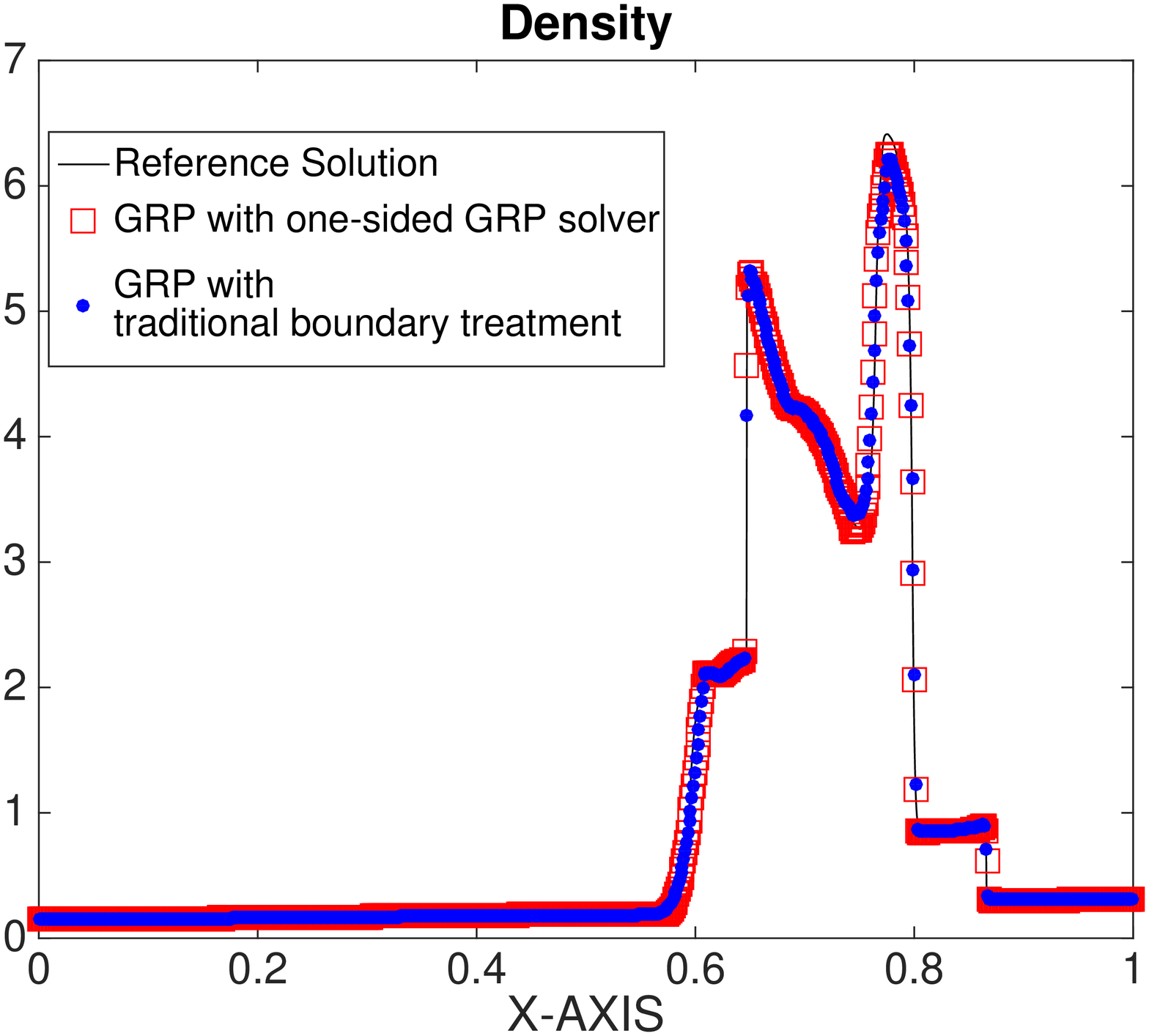}
}
\caption{\scriptsize The Woodward-Colella problem computed with the one-sided GRP solver (squares) and the traditional reflective boundary condition treatment (dots) with 400 cells (left) and 800 cells (right). The numerical scheme used in the interior domain is the GRP scheme. The solid lines are the reference solution computed with 4000 cells.}
\end{figure}
\vspace{0.2cm} 

\noindent
{\bf Example 4. The nozzle flow. } The nozzle flow problem is a classical quasi one-dimensional problem. Consider a flow in a converging-diverging nozzle  occupying the domain $x\in[0,1]$.  The cross-sectional area function $A(x)$ of the duct is given by
\begin{equation}
A(x)=\left\{
\begin{array}{ll}
\displaystyle A_{\rm in}{\rm exp}\left(-{\rm log}(A_{\rm in}){\rm sin}^{2}(2\pi x)\right),\quad 0\leq x\leq 0.25,\\[5pt]
\displaystyle A_{\rm ex}{\rm exp}\left(-{\rm log}(A_{\rm ex}){\rm sin}^{2}(\frac{2\pi (1-x)}{3})\right),\quad 0.25<x\leq 1
\end{array}
\right.
\end{equation}
with $A_{\rm in}=4.864317646$ and $A_{\rm ex}=4.234567901$. The governing equations  are the Euler equations with geometric source term \eqref{eq:Euler-1},\eqref{eq:euler-source}. Set $x=0$ as the entrance of the duct and $x=1$ as the exit. We are concerned with the present boundary treatment to attain the steady state solution.

Two types of steady states are discussed: A continuous steady state and a discontinuous steady state containing a standing shock wave. The initial data for both cases can take as
\begin{equation}
{\bf u}(0,x)=\left\{
\begin{array}{ll}
\displaystyle (\rho_0,0,p_0),\quad  x< 0.25,\\[5pt]
\displaystyle (\rho_0,0,\rho_0(p_{\rm ex}/p_0)^{\gamma}),\quad x>0.25,
\end{array}
\right.
\label{eq:duct-initial}
\end{equation}
where $\gamma=1.4$ and $\rho_0, p_0$ are parameters to be determined, $p_{\rm ex}$ is a constant value determined by the steady solution at $x=1$. In the previous study \cite{BenF,BenL1}, the inflow density, velocity and pressure are assigned to the inflow boundary condition, the outflow pressure is assigned as the outflow boundary condition. Here we apply the one-sided GRP solver to test its ability of attaining steady solutions.

For the first case, we set $\rho_{0}=p_{0}=1$ and $p_{\rm ex}=0.0272237$ in \eqref{eq:duct-initial}. This produces an isentropic continuous steady solutions which is defined by
\begin{equation}
\begin{array}{lll}
\displaystyle \rho(x)=\rho_0\left(1+\frac{\gamma-1}{2}M^{2}(x)\right)^{-\frac{1}{\gamma-1}},\\[8pt]
\displaystyle p(x)=p_0\left(1+\frac{\gamma-1}{2}M^{2}(x)\right)^{-\frac{\gamma}{\gamma-1}},\\[12pt]
\displaystyle v(x)=M(x)\sqrt{\gamma p(x)/\rho(x)},
\end{array}
\label{eq:duct-exact}
\end{equation}
in which the Mach number $M(x)=v(x)/c(x)$ is determined by $A(x)$ through the algebraic relation
\begin{equation}
\displaystyle A^{2}(x)=\frac{1}{M^{2}(x)}\left(\frac{2}{\gamma+2}\left(1+\frac{\gamma-1}{2}M^{2}(x)\right)\right)^{\frac{\gamma+1}{\gamma-1}}.
\end{equation}
In this case, the flow is transonic across the throat at the position $x=0.25$. Thus the inflow boundary condition at the entrance $x=0$ should be prescribed by
\begin{equation}
\begin{array}{ll}
\displaystyle p_{\rm in}:=p_0\left(1+\frac{\gamma-1}{2}M^{2}(0)\right)^{-\frac{\gamma}{\gamma-1}},\\
\displaystyle \rho_{\rm in}:=\rho_0\left(1+\frac{\gamma-1}{2}M^{2}(0)\right)^{-\frac{1}{\gamma-1}}.
\end{array}
\end{equation}
While at the exit $x=1$, the flow is supersonic and no boundary condition is needed. The computational result is given in Fig. 4.7 where 22 cells are used. The CFL number is 0.6 and the output time is $t=5$. The solution obtained by implementing the one-sided GRP solver converges to the exact steady one and is comparable with the result obtained in \cite{BenL1}. 

\begin{figure}[h]
\begin{minipage}[t]{\textwidth}
\centering
\includegraphics[width=\textwidth]{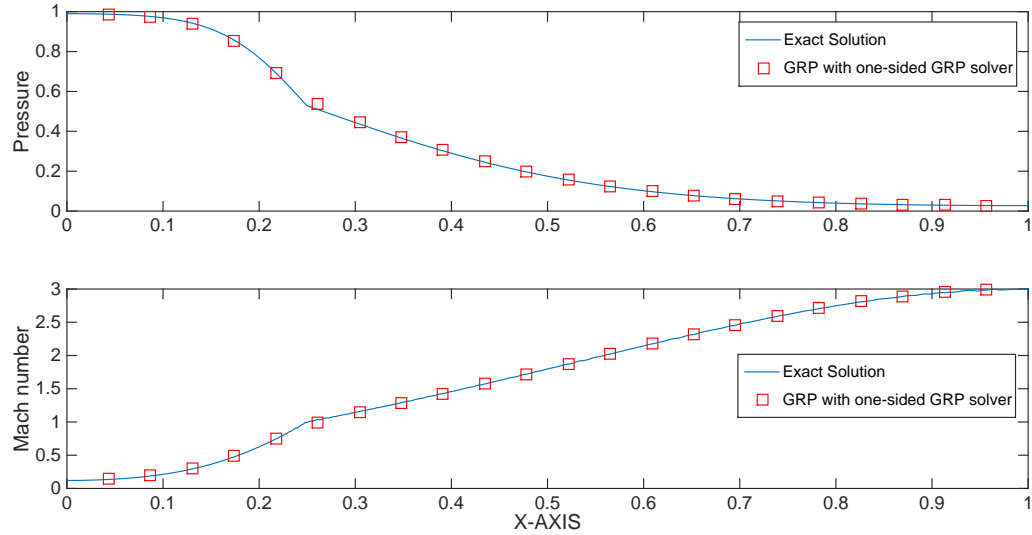}
\end{minipage}
\caption{ \scriptsize The computation of nozzle flow equations with continuous steady solutions by using the one-sided GRP solver. The pressure and Mach number at $t=5$ are shown with 22 cells. The solid line represents the exact solution given by \eqref{eq:duct-exact}.}
\end{figure}

\begin{figure}[h]
\begin{minipage}[t]{\textwidth}
\centering
\includegraphics[width=\textwidth]{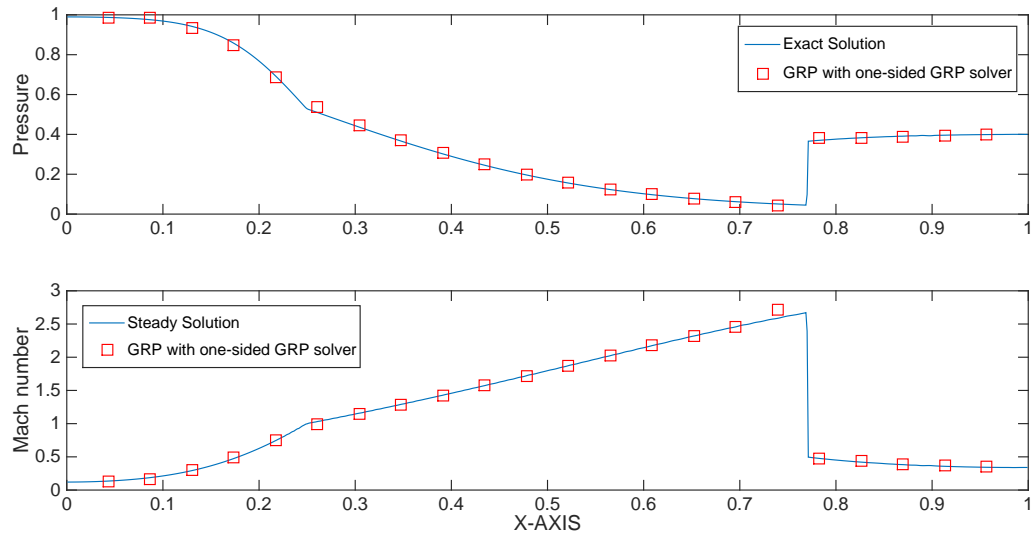}
\end{minipage}
\caption{\scriptsize The computation of nozzle flow equations with a standing shock wave by using the one-sided GRP solver. The pressure and Mach number at $t=5$ are shown with 22 cells. The solid line represents the exact solution given by \eqref{eq:duct-exact}.}
\end{figure}
\vspace{0.2cm} 

For the other case, where the steady solution contains a standing shock wave, we set $\rho_{0}=p_{0}=1$ and $p_{\rm ex}=0.4$ in \eqref{eq:duct-initial} to get the initial data. In this case, the flow jumps from supersonic to subsonic after passing the standing shock wave. As  the outflow is subsonic in this case, both inflow boundary condition and outflow boundary condition should be imposed. The inflow boundary condition is $\rho_0=p_0=1$ at the entrance $x=0$ and the outflow boundary condition is $p_{\rm ex}=0.4$ at the exit $x=1$. The computational result with 22 cells is given in Fig. 4.8. The CFL number is 0.6 and the output time is $t=5$. The solution obtained by taking the one-sided GRP solver matches well with the exact solution.

\begin{figure}[h]
\includegraphics[width=0.8\textwidth]{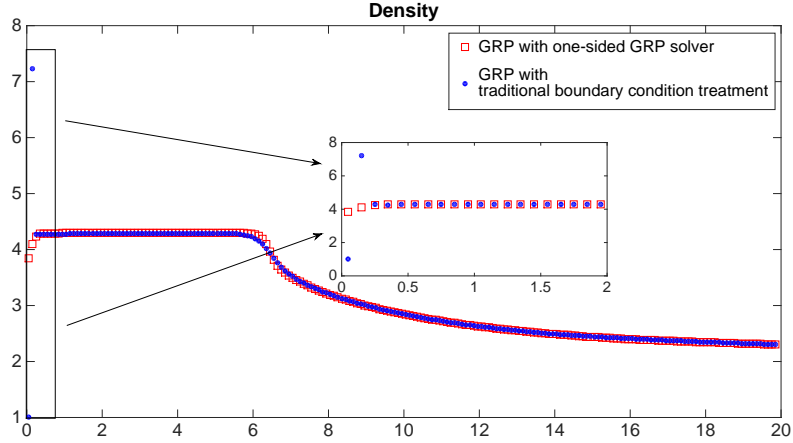}
\caption{ \scriptsize A spherical shock wave interacts with the symmetric center. We compute the density profile by the one-sided GRP solver (squares) and the traditional boundary condition treatment (dots). 200 cells are used.}
\end{figure}
\vspace{0.2cm} 

\noindent
{\bf Example 5. The spherical symmetric shock interaction problem.}  We test the one-sided GRP solver for the simulation of  the spherical symmetric flows where 
 a spherical shock wave interacts with the symmetric center. The initial data is taken to be 
 \begin{equation}
(\rho,v,p)(0,x)=\left\{
\begin{array}{ll}
(1.0,0.0,1.0/1.4),\quad 0\leq x\leq 2.0,\\[5pt]
(1.69997,-0.578906,1.528199),\quad 2.0<x\leq 10.0,
\end{array}
\right.
\end{equation}
such that a left-going spherical shock moves toward the center. The output time is $t=5.0$ with the CFL=0.5. One can see from Fig. 4.9 that near the symmetric center, the one-sided GRP solver has much better numerical performance compared with   the reflective boundary condition.  For more details about the GRP solver of radially symmetric flows, we refer to \cite{LLS} and references cited therein.

\n {\bf Example 6. Noh problem. } The Noh problem \cite{Noh} is a typical radially symmetric compressible flows problem. The governing equations include source term, which can be used to test the performance of the one-sided GRP solver. We consider the spherically converging flow of zero-pressure gas with $\gamma=5/3$. The initial data has the uniform form
\begin{equation}
[\rho,v,p]=[1,-1,0],\quad 0<r\leq 100,
\end{equation}
here $r$ is the radius. The exact solution consists of an expanding shock wave which begins from the center $r=0$. Here the initial pressure is set to be $10^{-6}$ instead of zero. The boundary condition at the rightmost cell is given by 
\begin{equation}
[\rho,v,p]^{n+1}(r)=[(1+t_{n+1}/r)^{2},-1,10^{-6}],\quad r\in[r_{K-1/2},r_{K+1/2}],
\end{equation}
which is the exact solution at $t=t_{n+1}$. On the left boundary one has $v(0,t)=0$. The one-sided GRP solver is implemented on both boundaries. The result is shown in Fig. 4.10. The discrepancies near the center is caused by the ``startup" of the captured shock wave, as pointed out in \cite{LLS}. The result obtained here has less
oscillations near the boundary compared with that in \cite{LLS}.
\vspace{2mm}

\begin{figure}[htbp]
\centering
\includegraphics[width=1.0\textwidth]{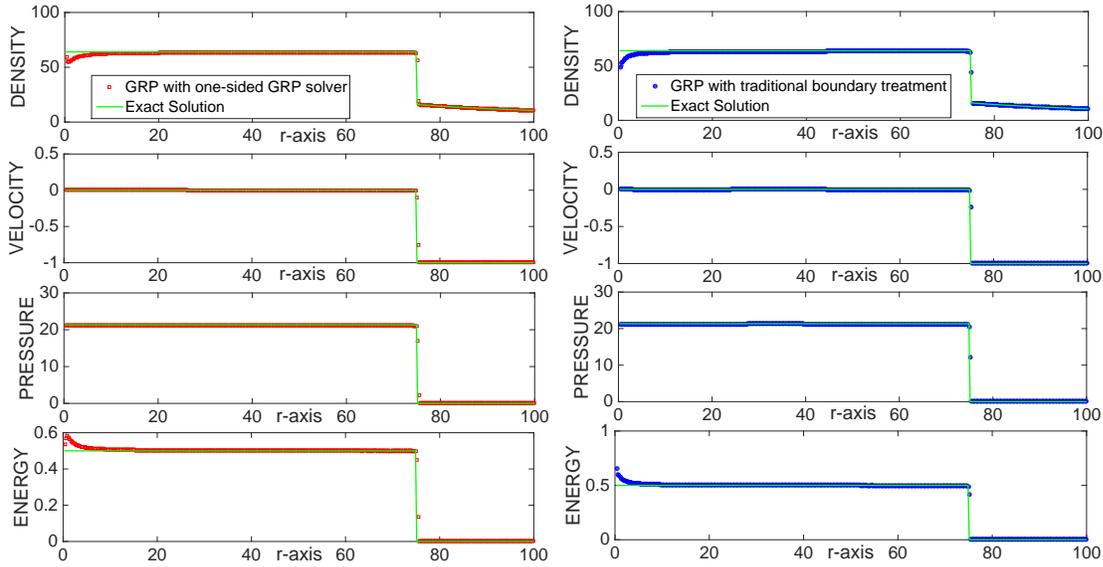}
\caption{\scriptsize The numerical results for Noh problem with 400 cells by using the one-sided GRP solver (squares) and the traditional reflective boundary treatment (dots), respectively, the solid line is given as the exact solution.}
\end{figure}

\n {\bf Example 7. The spherical explosion problem.} This is another problem of radially symmetric compressible flows. The initial gas is at rest with $\rho=21.7333, p=15.514$ for $0\leq r\leq 5$ and $\rho=2.0$, $p=1.0$ for $5\leq r\leq 50$. The spherical explosion is quite complex and a complete analysis can be found in \cite{LLS}. The numerical results are shown in Fig. 4.11, where we implement the GRP with two boundary condition treatments : the one-sided GRP solver and the method  developed in \cite{LLS}. From Fig. 4.11 one see that the one-sided GRP solver has a good agreement with the method that proposed in \cite{LLS}.
\vspace{2mm}

\begin{figure}[htbp]
\begin{minipage}[t]{\textwidth}
\centering
\includegraphics[width=0.95\textwidth]{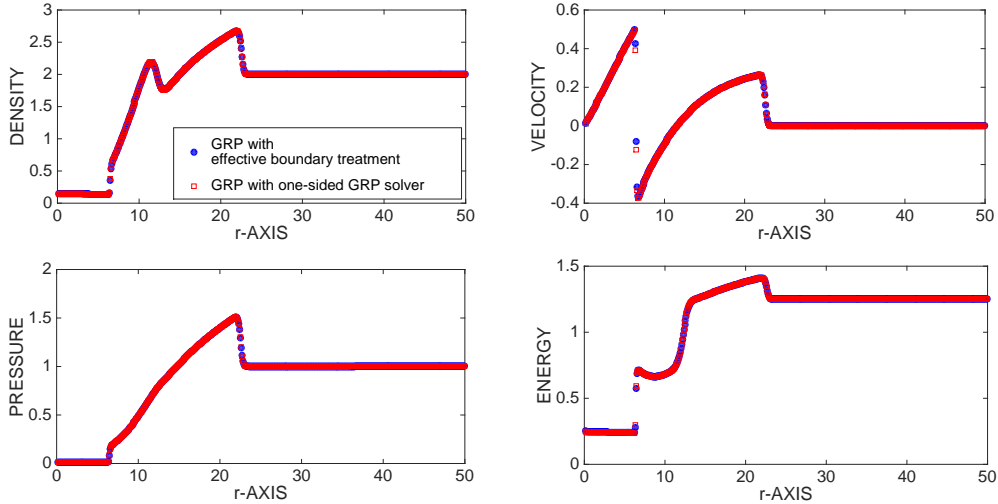}
\end{minipage}
\caption{\scriptsize The comparison of the results of the spherical explosion problem with two boundary condition treatments: one-sided GRP solver (squares) and the effective boundary condition treatment (dots) proposed in \cite{LLS}.}
\end{figure}

\begin{figure}[htbp]
\begin{minipage}[t]{\textwidth}
\centering
\includegraphics[width=0.975\textwidth]{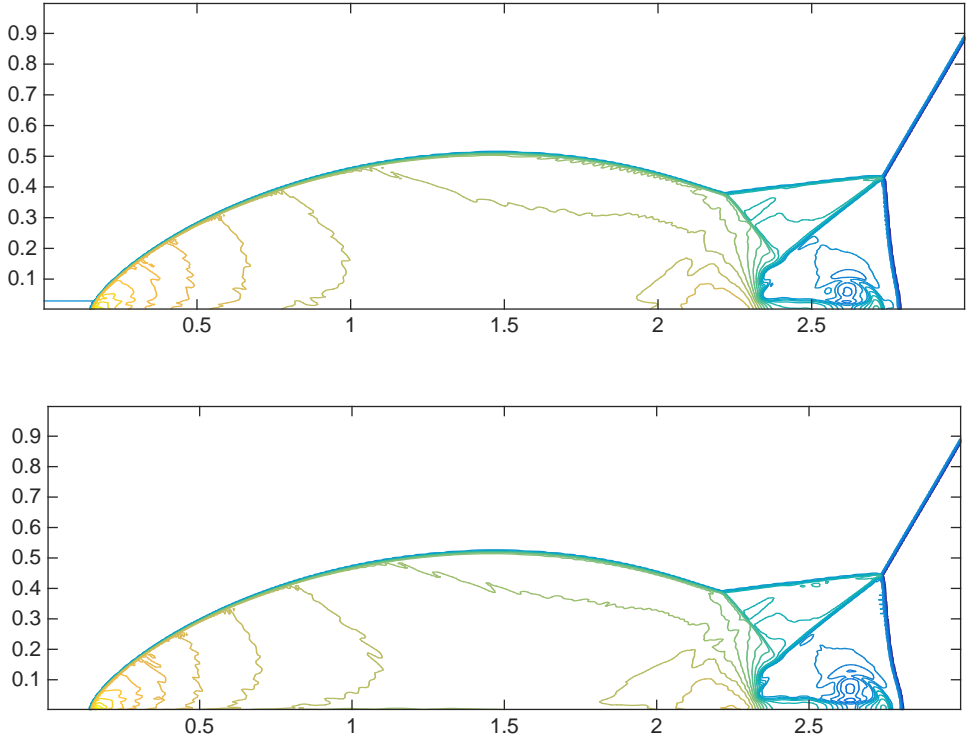}
\end{minipage}
\caption{\scriptsize The numerical results of the double Mach reflection problem. The upper is the GRP scheme with the one-sided GRP solver. The lower is the GRP scheme with the traditional reflective boundary condition treatment.}
\end{figure}

\n {\bf Example 8. The double Mach reflection problem.} We turn to two-dimensional example. The computational domain is $[0,4] \times [0,1]$, and $[0,3] \times [0,1]$ is shown. A solid-wall is at the bottom of the domain starting from $x=\frac16$. Initially a Mach $10$ shock wave is moving to the right which is at the position $x=\frac16, y=0$ and makes $\frac{\pi}{3}$ angle with the $x$-axis. More details about the problem can be seen in \cite{WC}.

We compute the problem by using the traditional boundary condition treatment and the one-sided GRP solver, respectively, to deal with the reflective boundary condition along the bottom wall $\{(x,y):\frac16<x<4,y=0\}$. The results are displayed in Fig. 4.12 with $30$ contours of the density at time $t=0.2$ where $720\times180$ cells are used here. The CFL number is $0.6$. From the figure we see that the one-sided GRP solver works well for the two-dimensional solid-wall boundary condition.

\begin{figure}[htbp]
\begin{minipage}[t]{\textwidth}
\centering
\includegraphics[width=0.975\textwidth]{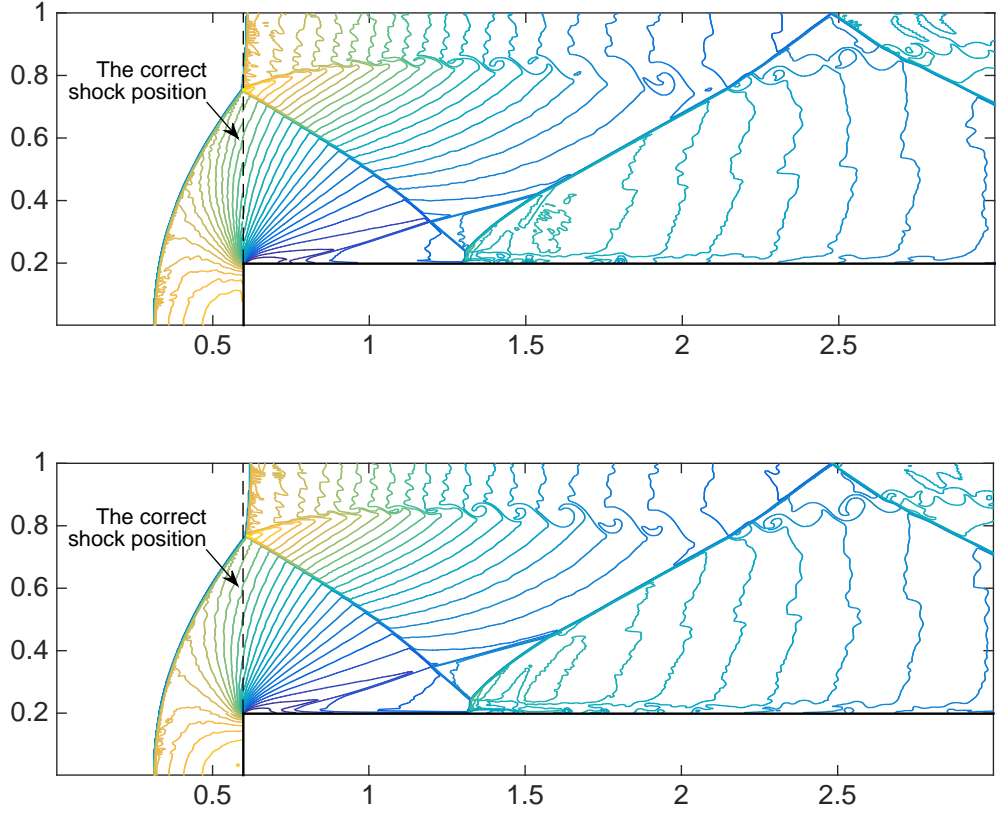}
\end{minipage}
\caption{\scriptsize The numerical results of the forward facing step problem. The upper is the GRP scheme with the one-sided GRP solver. The lower is the GRP scheme with the traditional reflective boundary condition treatment.}
\end{figure}

{\bf Example 9. The forward facing step problem.} This is another classical test problem for the two-dimensional equations. The wind tunnel is $1$ length unit wide and $3$ length units long. The step is $0.2$ length units high and is located $0.6$ length units from the left-hand end of the tunnel. Initially a unit right moving Mach 3 shock wave with $(\rho_0,v^{x}_0,v^{y}_0,p_0)=(1.4,3,0,1)$ in the tunnel. The reflective boundary conditions are applied along all the walls.

Again, we compute the problem by using the traditional boundary condition treatment and the one-sided GRP solver, respectively. The CFL number is $0.6$. The results are displayed in Fig. 4.13 with 900$\times$300 cells at time $t=4$. A three-shock Mach reflection configuration is formed. According to \cite{WC}, the correct Mach stem is located at $x=0.6$.  We can see that the results obtained by the one-sided GRP solver has the shock at the correct position, compared with  that obtained by the reflective boundary condition treatment.

\appendix

\section{The useful one-sided GRP coefficients.} The coefficients of the one-sided GRP solver are collected in Table 1. In this table, the $1$-shock ($3$-shock, resp.) refers to the shock associated with the $v-c$ characteristic family ($v+c$, resp.). The same for the $1$-rarefaction wave and $3$-rarefaction wave. We deal with both the left boundary case and the right boundary case. For the left boundary case, the one-sided Riemann problem \eqref{eq:RP} has the solution which consists of a single $3$-shock wave or a $3$-rarefaction wave. Similarly, when there exists a right boundary, the solution of \eqref{eq:RP} consists of a single $1$-shock wave or a $1$-rarefaction wave. 

\begin{table}[]
\scriptsize
\centering
\caption{The coefficients in \eqref{3.7} for all possible cases.}
\label{Tab03}
\begin{tabular}{cc} 
\toprule
 1-rarefaction wave    \qquad  \qquad  \qquad    \qquad   \qquad   \qquad  & $(a_{L},b_{L})=(a_{L}^{\rm rare},b_{L}^{\rm rare}),d_{L}=d_{L}^{\rm rare}$          \\ [6pt]
1-shock wave       \qquad   \qquad  \qquad   \qquad  \qquad    \qquad    &$(a_{L},b_{L})=(a_{L}^{\rm shock},b_{L}^{\rm shock}),d_{L}=d_{L}^{\rm shock}$          \\  [6pt]
3-rarefaction wave        \qquad \qquad   \qquad   \qquad   \qquad   \qquad    &$(a_{R},b_{R})=(a_{R}^{\rm rare},b_{R}^{\rm rare}),d_{R}=d_{R}^{\rm rare}$         \\  [6pt]
3-shock wave      \qquad    \qquad  \qquad   \qquad  \qquad    \qquad  &$(a_{R},b_{R})=(a_{R}^{\rm shock},b_{R}^{\rm shock}),d_{R}=d_{R}^{\rm shock}$    \\  
\bottomrule
\end{tabular}
\end{table}

We denote $J=L$ or $R$ in the rest part of the paper. The one-sided Riemann solution is denoted as ${\bf u}_{J}^{*}$ which can be obtained through solving \eqref{eq:RP}.
The Riemann invariants $\phi$ and $\psi$ are introduced by
\begin{equation}\label{A.2}
\psi = v+\frac{2c}{\gamma-1},\quad \phi = v-\frac{2c}{\gamma-1}.
\end{equation}
For the second law of thermodynamics, we have 
\begin{equation}
TdS=\frac{{\rm d}p}{(\gamma-1)\rho}-\frac{c^2}{(\gamma-1)\rho}{\rm d}\rho.
\end{equation}
These are useful to derive the coefficients below. More details can be found in \cite{BenL1}.

{\bf A.1.} (Nonsonic case) The coefficients for rarefaction waves are given by
\begin{equation}\label{A.3}
\begin{array}{lll}
\displaystyle (a_{ L}^{\rm rare},b_{ L}^{\rm rare})=\left(1,-\frac{1}{\rho_{ L}^{*}c_{ L}^{*}}\right),\quad (a_{ R}^{\rm rare},b_{ R}^{\rm rare})=\left(1,-\frac{1}{\rho_{ R}^{*}c_{ R}^{*}}\right),\\
 d_{J}^{\rm rare} =  \left[\displaystyle \frac{1+\mu^2}{1+2\mu^2}\theta_{J}^{\frac{1}{2\mu^2}}+\frac{\mu^2}{1+2\mu^2}\theta_{J}^{\frac{1+\mu^2}{\mu^2}}\right]T_{J}S'_{J} + {sgn}(J)c_{J}\left(\eta'(J) + \frac{a'(0)}{a(0)}v_{J}\right)\theta_{J}^{\frac{1}{2\mu^2}} \\[12pt]
 \quad \quad    + \displaystyle  \frac{a'(0)}{a(0)}c_{J}^{*}(\Phi_{J}+sgn(J)v_{J}^{*}),
\end{array}
\end{equation}
where $\displaystyle \mu^2=\frac{\gamma-1}{\gamma+1}$, $ \displaystyle \theta_{ L}=\frac{c_{ L}^{*}}{c_{ L}}$, $\displaystyle \theta_{ R}=\frac{c_{R}^{*}}{c_{R}}$, and $\Phi_{J}$ are given by
\begin{equation}
\Phi_{J}=\left\{
\begin{array}{lll}
\displaystyle \frac{(\mu^2-1)c_{J}^{*}}{\mu^{2}(4\mu^2-1)}\left[1-\theta_{J}^{\frac{1-4\mu^2}{2\mu^2}}\right]-sgn(J)\frac{\eta(J)}{2\mu^2-1}\left[1-\theta_{J}^{\frac{1-2\mu^2}{2\mu^2}}\right], {\rm if~} \gamma \neq \frac53, \gamma \neq 3,\\[10pt]
\displaystyle c_{J}-c_{J}^{*}-sgn(J)\eta(J){\rm log}(\theta_{J}), \quad {\rm if}~ \gamma=3,\\
\displaystyle -2[3c_{J}^{*}{\rm log}\theta_{J}-sgn(J)\eta(J)(1-\theta_{J})],\quad {\rm if~\gamma=\frac53}.
\end{array}
\right.
\end{equation}
Here
\begin{equation}\label{A.5}
sgn(J)=\left\{\begin{array}{ll}
-1,\quad &{\rm if~}J=L,\\
1,\quad &{\rm if~}J=R,
\end{array}
\right.
\quad
\eta(J)=\left\{\begin{array}{ll}
\psi_{L},\quad &{\rm if~}J=L,\\
\phi_{R},\quad &{\rm if~}J=R.
\end{array}
\right.
\end{equation}

The coefficients for shock waves are given by
\begin{equation}
\begin{array}{ll}
\displaystyle  a_{J}^{\rm shock}=1+sgn(J)\rho_{J}^{*}(\sigma_{J}-v_{J}^{*})\Phi_{1}^{J},\quad b_{J}^{\rm shock}=-\frac{1}{\rho_{J}^{*}(c_{J}^{*})^{2}}(\sigma_{J}-v_{J}^{*})-sgn(J)\Phi_{1}^{J},\\
\displaystyle  d_{J}^{\rm shock}=L_{\rho}^{J}\rho'_{J}+L_{p}^{J}p'_{J}+L_{v}^{J}v'_{J}-\frac{a'(0)}{a(0)}j_{R}.
\end{array}
\end{equation}
All the variables involved are 
\begin{equation}
\begin{array}{lll}
\displaystyle \sigma_{J}=\frac{\rho_{J}^{*}v_{J}^{*}-\rho_{J}v_{J}}{\rho_{J}^{*}-\rho_{J}},\\[9pt]
\displaystyle L_{\rho}^{J}=sgn(J)(\sigma_{J}-v_{J})\Phi_{3}^{J},\quad  L_{v}^{J}=\sigma_{J}-v_{J}-sgn(J)(\rho_{J}c_{J}^{2}\Phi_{2}^{J}+\rho_{J}\Phi_{3}^{J}),\\
\displaystyle L_{p}^{J}= -\frac{1}{\rho_{J}}+sgn(J)(\sigma_{J}-v_{J})\Phi_2^{J},\quad j_{R}=sgn(J)\rho_{J}v_{J}(c_{J}^{2}\Phi_{2}^{J}+\Phi_{3}^{J})-(\sigma_{J}-v_{J}^{*})v_{J}^{*}.
\end{array}
\end{equation}
Here $H_{i}^{J}=H_{i}(p_{*}^{J};p_{J},\rho_{J}), i=1,2,3$, $H_{i}$ is given by
\begin{equation}
\begin{array}{lll}
\displaystyle H_1(p;\bar{p},\bar{\rho})=\frac12\sqrt{\frac{1-\mu^2}{\bar{\rho}(p+\mu^2\bar{p})}}\frac{p+(1+2\mu^2)\bar{p}}{p+\mu^2\bar{p}},\\[7pt]
\displaystyle H_2(p;\bar{p},\bar{\rho})=-\frac12\sqrt{\frac{1-\mu^2}{\bar{\rho}(p+\mu^2\bar{p})}}\frac{(2+\mu^2)p+\mu^{2}\bar{p}}{p+\mu^2\bar{p}},~\displaystyle H_3(p;\bar{p},\bar{\rho})=-\frac{p-\bar{p}}{2\bar{\rho}}\sqrt{\frac{1-\mu^2}{\bar{\rho}(p+\mu^2\bar{p}}}.
\end{array}
\end{equation}

Denote $D/Dt = \partial/\partial t + v \partial /\partial x$. Then we have 
\begin{equation}\label{A.9}
\begin{array}{ll}
\displaystyle \frac{{\rm \partial} v}{{\rm \partial} t}=\frac{{\rm D}v}{{\rm D}t}+\frac{v}{\rho c^2}\frac{{\rm D}p}{{\rm D}t}+\frac{a'(0)}{a(0)}v^2,\\[12pt]
\displaystyle \frac{{\rm \partial}p}{{\rm \partial}t}=\frac{{\rm D}p}{{\rm D}t}+\rho v\frac{{\rm D}v}{{\rm D}t}.
\end{array}
\end{equation}
Remember that we always have $(\partial v/\partial t)^{*}=g'(t)$ on the boundary, then the expected instantaneous values $(\partial {\bf u}/\partial t)^{*}$ can be obtained directly through solving \eqref{3.7} with \eqref{A.9}

{\bf A.2.} (Sonic case) When the left boundary is located inside the 3-rarefaction wave, we have 
\begin{equation}\label{A.10}
\left(\frac{{\rm\partial} v}{{\rm \partial}t}\right)^{*}=g'(t),\quad \left(\frac{\partial p}{\partial t}\right)^{*}=\rho_{R}^{*}v_{R}^{*}\left[\left(\frac{{\rm\partial} v}{{\rm \partial}t}\right)^{*}-\theta^{\frac{2\gamma}{\gamma-1}}T_{R}S'_{R}-\frac{a'(0)}{a(0)}(v_{R}^{*})^{2}\right].
\end{equation}
 Similarly, when the right boundary locates inside the 1-rarefaction wave, we just replace ${\bf u}_{R}^{*}$, $T_{R}S'_{R}$ by ${\bf u}_{L}^{*}$, $T_{L}S'_{L}$ in \eqref{A.10}. 
 
{\bf A.3.} (Acoustic case) Assume that  on the right boundary, ${\bf u}_{L}^{*}={\bf u}_{L}, ({\bf u}_{L}^{*})'\neq {\bf u}'_{L}$, or on the left boundary, ${\bf u}_{R}^{*}={\bf u}_{R}, ({\bf u}^{*}_R)'\neq {\bf u}'_{R}$, we have the acoustic case. $(\partial v/\partial t)^{*}$ and $(\partial p/\partial t)^{*}$ can be given by
\begin{equation}
\begin{array}{ll}
\displaystyle \left(\frac{{\rm\partial} v}{{\rm \partial}t}\right)^{*}=g'(t),\\
 \displaystyle  \left(\frac{\partial p}{\partial t}\right)^{*}=sgn(J)\rho_{J}^{*}c_{J}^{*}-\rho_{J}^{*}\left(v_{J}^{*}-sgn(J)c_{J}^{*}\right)\left[\frac{p_{J}^{'}}{\rho_{J}^{*}}-sgn(J)c_{J}^{*}v'_{J}\right]-\frac{a'(0)}{a(0)}(\rho_{J}^{*}v_{J}^{*})^{3}.
\end{array}
\end{equation}
And the quantity $(\partial \rho/\partial t)^{*}$ is calculated from the EOS,
\begin{equation}
\displaystyle \left(\frac{\partial \rho}{\partial t}\right)^{*}=\frac{1}{(c_{J}^{*})^2}\left[\left(\frac{\partial p}{\partial t}\right)^{*}+v_{J}^{*}\left(p'_{J}-(c_{J}^{*})^{2}\rho'_{J}\right)\right],
\end{equation}
where $J$ takes $L$ or $R$, the definition of $sgn(J)$ is referred to \eqref{A.5}.

\section*{Acknowledgment} The first author is supported by NSFC (Nos. 11771054, 12072042,91852207), the Sino-German Research Group Project (No. GZ1465)  and Foundation of LCP.

\bibliographystyle{amsplain}

\begin{thebibliography}{20}

\bibitem{Bank-Ben-Artzi}
M. Bank and M.  Ben-Artzi,  Scalar conservation laws on a half-line: a parabolic approach. {\em J. Hyper. Diff. Equat. } 7 (2010) 165--189. 

\bibitem{Bardos-1979} C. Bardos, A. Y. le Roux, and J.-C.  N\'ed\'elec,  First order quasilinear equations with boundary conditions. {\em Comm. Partial Diff. Equat.} 4 (1979) no. 9, 1017--1034.



\bibitem{BenF} M. Ben-Artzi and  J. Falcovitz, Generalized Riemann Problems in Computational Fluid Dynamics, {\em Cambridge University Press}, 2003.

\bibitem{BenL1} M. Ben-Artzi and  J. Li, Hyperbolic balance laws: Riemann invariants and the generalized Riemann problem, {\em Numer. Math.}, 106 (2007) 369-425.

\bibitem{BenLW} M. Ben-Artzi, J. Li and  G. Warnecke,  A direct Eulerian GRP scheme for compressible fluid flows, {\em J. Comput. Phys., }218 (2006) 19-43.

\bibitem{Ben-Dor-2007} G. Ben-Dor, Shock wave reflection phenomena, {\em Springer}, 2007.


\bibitem{Berger} M. Berger, C. Helzel and  R. LeVeque, h-Box methods for the approximation of hyperbolic conservation laws on irregular grids, {\em  SIAM J. Numer. Anal.}, 41 (2003) 893--918.

\bibitem{CF} R. Courant and K.O. Friedrichs, Supersonic flow and shock waves, {\em Springer, New York}, 1948.

\bibitem{CH} T. Chang and L. Hsiao, The Riemann problem and interaction of waves in gas dynamics. {\em Pitman Monographs and Surveys in Pure and Applied Mathematics, 41. Longman Scientific and Technical, Harlow}, 1989. 

\bibitem{Dafermos-book} C. Dafermos, Hyperbolic Conservation Laws in Continuum Physics, {\em volume 325 of A series of comprehensive studies in Mathematics, Springer-Verlag, Berlin,} 2000. 


\bibitem{DDJ} G. Dakin, B. Despr\'{e}s and S. Jaouen, Inverse Lax-Wendroff boundary treatment for compressible Lagrange-remap hydrodynamics on Cartesian grids, {\em J. Comput. Phys.}, 353 (2018) 228-257.


\bibitem{DL1}  Z. Du and  J. Li,  A two-stage fourth order time-accurate discretization for Lax-Wendroff type flow solvers {\rm I\!I}. High order numerical boundary conditions, {\em J. Comput. Phys.,} 369 (2018) 125-147.


\bibitem{DL3}  Z.  Du and J. Li,  Accelerated Piston Problem and High Order Moving Boundary Tracking Method for Compressible Fluid Flows, {\em SIAM J. Sci. Comput.}, A1558-A1581.

\bibitem{FAH}  J. Falcovitz, G. Alfandary and G. Hanoch, A two-dimensional conservation laws scheme for compressible flow with moving boundaries, {\em J. Comput. Phys.}, 138 (1997) 83-102.

\bibitem{FJ}  H. Forrer and R. Jeltsch, A high-order boundary treatment for Cartesian-grid methods, {\em J. Comput. Phys.}, 140 (1998) 259-277.

\bibitem{Kreiss-book} H.-O. Kreiss and  J. Lorenz, Initial-Boundary Value Problems and the Navier--Stokes Equations, {\em Academic Press, San Diego}, 1989.

\bibitem{KP1}  H.-O. Kreiss and N. Petersson, A second order accurate embedded boundary method for the wave equation with Dirichlet data, {\em SIAM J. Sci. Comput.}, 27 (2006) 1141-1167.

\bibitem{KP2}  H.-O. Kreiss, N. Petersson and J. Ystr\"{o}m, Difference approximations of the Neumann problem for the second order wave equation, {\em SIAM J. Numer. Anal.}, 42 (2004) 1292-1323.

\bibitem{KB}  L. Krivodonova and M. Berger, High-order accurate implementation of solid wall boundary conditions in curved geometries, {\em J. Comput. Phys.}, 211 (2006) 492-512.

\bibitem{Lax-1957} P. Lax, Hyperbolic systems of conservation laws II, {\em Comm. Pure Appl. Math.}, 10 (1957)  537-566. 

\bibitem{LL} X. Lei and J. Li, Transversal effects of high order numerical schemes for compressible fluid flows, {\em App. Math. Mech. Engl. Ed.}, 40 (2019) 343-354.

\bibitem{LLS}  J. Li, T. Liu and Z. Sun, Implementation of the GRP scheme for computing radially symmetric compressible fluid flows, {\em J. Comput. Phys.}, 228 (2009) 5867-5887.

\bibitem{Li1}  J. Li,  Two-stage fourth order: Temporal-spatial coupling in computational fluid dynamics (CFD), {\em Adv. Aerodynam.}, 1 (2019) 1:3; 1-36.

\bibitem{LD}  J. Li and Z. Du,  A two-stage fourth order temporal discretization for on the Lax-Wendroff type flow solvers {\rm I}. Hyperbolic conservation laws, {\em SIAM, J. Sci. Comput.}, 38 (2016) 3046-3069.

\bibitem{Li-Yu} Ta-Tsien Li and Wen-Ci Yu,  Boundary value problems for quasilinear hyperbolic systems. {\em Duke University Mathematics Series, V. Duke University, Mathematics Department, Durham, NC,} 1985.

\bibitem{Noh}  W. F. Noh, Errors for calculations of strong shocks using an artificial viscosity and artificial heat flux, {\em J. Comput. Phys.}, 72 (1987) 78-120.


\bibitem{P-L-1992} T. J. Poinsot and S. K. Lee, Boundary conditions for direct simulations of compressible viscous flows, {\em J. Comput. Phys.,} 101 (1992) 102-129. 

\bibitem{SP} B. Sj\"ogreen and N. Petersson, A Cartesian embedded boundary method for hyperbolic conservation laws, {\em Commun. Comput. Phys.}, 2 (2007) 1199-1219.

\bibitem{TS1} S. Tan and C.-W. Shu, Inverse Lax-Wendroff procedure for numerical boundary conditions of conservation laws, {\em J. Comput. Phys.}, 229 (2010) 8144-8166.

\bibitem{TWSN} S. Tan, C. Wang, C.-W. Shu and J. Ning, Efficient implementation of high order inverse Lax-Wendroff boundary treatment for conservation laws, {\em J. Comput. Phys.,} 231 (2012) 2510-2527.

\bibitem{Toro-book} E. Toro, Riemann solvers and numerical methods for fluid dynamics. A practical introduction. Third edition. {\em Springer-Verlag, Berlin}, 2009. 

\bibitem{WC} P. Woodward and P. Collela, The numerical simulation of two-dimensional fluid flow with strong shocks, {\em J. Comput. Phys., } 54(1984) 115-173.

\end{thebibliography}

\end{document}